\documentclass[11pt,leqno]{amsart}
\usepackage[T1]{fontenc}
\usepackage{lmodern,amsmath,amssymb,mathtools,mathrsfs,booktabs,array,longtable}
\usepackage{microtype,color}
\definecolor{DOIBlue}{rgb}{0,0.30588235294,0.6}
\usepackage{etoolbox}
\usepackage[margin=1.05in]{geometry}
\usepackage[colorlinks=true,linkcolor=blue,citecolor=red,urlcolor=DOIBlue,linktoc=all]{hyperref}
\usepackage[noabbrev]{cleveref}
\DeclareFontFamily{U}{msb}{}
\DeclareFontShape{U}{msb}{m}{n}{<->msbm10}{}
\DeclareFontFamily{U}{msa}{}
\DeclareFontShape{U}{msa}{m}{n}{<->msam10}{}

\makeatletter
\patchcmd{\@maketitle}{\dimen@34\p@}{\dimen@22\p@}{}{}
\patchcmd{\@starttoc}{32\p@\@plus14\p@}{14\p@\@plus3\p@}{}{}
\AtBeginDocument{\def\@textbottom{\vskip\z@\@plus\baselineskip}}
\patchcmd{\thebibliography}{\labelsep .5em}{\labelsep .8em}{}{}
\def\subsubsection{\@startsection{subsubsection}{3}%
  \z@{.5\linespacing\@plus.7\linespacing}{-.5em}{\normalfont\bfseries}}
\def\paragraph{\@startsection{paragraph}{4}%
  \z@\z@{-\fontdimen2\font}{\normalfont\bfseries}}
\newcommand{\introproofheading}{\@startsection{paragraph}{4}%
  \z@{.65\linespacing\@plus.15\linespacing\@minus.1\linespacing}%
  {.22\linespacing}{\normalfont\bfseries}}
\makeatother
\numberwithin{equation}{section}
\newtheorem{maintheorem}{Theorem}

\newtheorem{theorem}{Theorem}[section]
\newtheorem{lemma}[theorem]{Lemma}
\newtheorem{proposition}[theorem]{Proposition}
\newtheorem{corollary}[theorem]{Corollary}

\newtheoremstyle{bolddefinition}
  {\topsep}{\topsep}{\normalfont}{}{\bfseries}{.}{.5em}{}
\theoremstyle{bolddefinition}

\newtheorem*{remark*}{Remark}
\crefname{maintheorem}{Theorem}{Theorems}
\Crefname{maintheorem}{Theorem}{Theorems}
\crefname{definition}{Definition}{Definitions}
\crefname{lemma}{Lemma}{Lemmas}
\crefname{corollary}{Corollary}{Corollaries}

\newcommand{\C}{\mathbb C}
\newcommand{\Disc}{\mathbb D}
\newcommand{\ddc}{dd^c}
\DeclareMathOperator{\Ric}{Ric}
\DeclareMathOperator{\av}{av}
\DeclareMathOperator{\tr}{tr}
\DeclareMathOperator{\dis}{dist}

\newcommand{\R}{\mathbb R}
\newcommand{\T}{\mathbb T}
\DeclareMathOperator{\vol}{vol}
\allowdisplaybreaks[2]

\title[Brody hyperbolicity, canonical ampleness, and general type]
{Brody Hyperbolicity, Canonical Ampleness, and Subvarieties of General Type}
\author[X. Qin]{Xiangsen Qin}
\address{Chern Institute of Mathematics and LPMC, Nankai University,
Tianjin 300071, China}
\email{qinxiangsen@nankai.edu.cn}
\author[Y. Tian]{Yuanhong Tian}
\address{School of Mathematics and Statistics, Beijing Institute of Technology,
Beijing 100081, China}
\email{yuanhongtian@bit.edu.cn}
\hypersetup{pdftitle={Brody Hyperbolicity, Canonical Ampleness, and Subvarieties of General Type},
  pdfauthor={Xiangsen Qin and Yuanhong Tian}}

\begin{document}

\begin{abstract}
Let \(X\) be a compact Brody-hyperbolic complex manifold of positive dimension. We prove that its canonical bundle \(K_X\) is ample. It follows that \(X\) is projective, which answers Kobayashi's question on the projectivity of compact hyperbolic manifolds. We also consider a generically finite holomorphic map \(g\colon Z\to W\) from a positive-dimensional reduced irreducible compact complex space \(Z\) in Fujiki's class \(\mathcal C\) to a compact hyperbolic complex space \(W\). Then \(Z\) and \(g(Z)\) are of general type, where $g(Z)$ is equipped with its reduced structure. Together these results imply that every positive-dimensional irreducible closed analytic subspace of \(X\) is of general type. For projective varieties, including singular ones, the result for generically finite maps gives the forward implication of Lang's geometric conjecture. As an application, we show that every Hermitian metric of constant negative Chern holomorphic sectional curvature on a compact complex manifold is Kähler. 
\end{abstract}

\maketitle
\setcounter{tocdepth}{1}
\hypersetup{bookmarksdepth=3}
\tableofcontents

\section{Introduction}\label{sec:introduction}

Let $X$ be a complex space. A holomorphic map $\C\to X$ is called
an \emph{entire curve}. We call $X$ \emph{Brody hyperbolic} if every
entire curve in $X$ is constant. Throughout, \emph{hyperbolic} means
Brody hyperbolic.
For compact complex manifolds, this condition is equivalent to
Kobayashi hyperbolicity (the Kobayashi pseudodistance separates distinct points) by Brody's theorem \cite[Theorem~4.1]{Brody78}.
The equivalence also holds on each connected component of a compact
complex space \cite[Ch.~III, Theorem~2.1]{Lang87}. Unless otherwise
stated, all ambient complex spaces have positive dimension.

Let $S$ be a compact Riemann surface of genus $g$. By the
uniformization theorem \cite{Koebe07,Poin08}, $S$ is hyperbolic
if and only if $g\ge2$. Since $\deg K_S=2g-2$,
$K_S$ is ample if and only if $g\ge2$. By the Riemann--Roch theorem
\cite{Rie57,Roc65}, this is also equivalent to $S$ being of general type. Thus
\begin{equation}\label{eq:curve-hyperbolicity}
 S\text{ is hyperbolic}
 \quad\Longleftrightarrow\quad K_S\text{ is ample}
 \quad\Longleftrightarrow\quad S\text{ is of general type}.
\end{equation}

In 1970, Kobayashi asked whether every compact $C$-hyperbolic manifold,
or more generally every compact hyperbolic manifold, is projective
\cite[Ch.~IX, \S3]{Kob70} (see \cite[Problem~7, p.~93]{Kob73} for the Russian
translation). A manifold is $C$-hyperbolic if it
admits an unramified covering on which bounded holomorphic functions
separate points \cite[Ch.~IX, \S3]{Kob70}. By Kodaira's embedding
theorem \cite[Main Theorem, p.~314]{Kod54}, a compact complex manifold
$X$ is projective if it carries a positive holomorphic line bundle. This leads to
the stronger question of whether hyperbolicity forces the canonical
bundle $K_X$ itself to be ample.

In 1972, Kobayashi observed that he knew of no compact hyperbolic
manifold $X$ for which $c_1(X)$ fails to be negative
\cite[Ch.~III, \S2, pp.~83--84 and~87]{Kob72}. Here $c_1(X)<0$ means that the first Bott--Chern class of $T_X$ has a
negative definite representative. Since this is equivalent to ampleness of $K_X$, this observation concerns
precisely the ampleness question above. In 1986, Lang stated the smooth
projective case \cite[Conjecture~7.2]{Lang86}. Peternell stated Kobayashi's conjecture for smooth projective manifolds
and obtained partial results in dimension three
\cite[\S2, (2.2)--(2.3), pp.~316--317]{Pet91}. Campana attributed the question to Kobayashi for compact
K\"ahler manifolds in his 1991 conference paper \cite[p.~65]{Cam91}
and for compact connected complex manifolds in his 1992 paper
\cite[\S0, p.~1]{Cam92}. The converse fails
in higher dimensions. For example, the Fermat quintic surface in
$\mathbb P^3$ has ample canonical bundle but contains a line.

A compact manifold with ample canonical bundle is of general type. For a
singular space, general type can be tested on a resolution. In 1975,
Kobayashi and Ochiai suggested that every reduced irreducible compact
hyperbolic complex space $Z$ is of general type, meaning that
$\widetilde Z$ is of general type for any resolution $\widetilde Z\to Z$
\cite[\S1, p.~7 and footnote~1]{KO75}.

Every closed subvariety of a hyperbolic projective variety inherits
hyperbolicity. Lang conjectured that a possibly singular, irreducible projective
variety $X$ is hyperbolic if and only if every
positive-dimensional irreducible closed subvariety of $X$, including
$X$ itself, is of general type \cite[Conjecture~5.6]{Lang86}.

We prove that every positive-dimensional compact hyperbolic complex
manifold has ample canonical bundle, as suggested by Kobayashi's
1972 observation. Kodaira's embedding theorem then gives an affirmative
answer to his 1970 question about projectivity.
We also prove the general-type prediction of Kobayashi and Ochiai for
positive-dimensional reduced irreducible compact hyperbolic complex
spaces in Fujiki's class $\mathcal C$ (\S\ref{sec:comparison}).
Applied to closed subvarieties of hyperbolic projective varieties, this
gives the forward implication of Lang's conjecture.

\subsection{Canonical ampleness}

\begin{maintheorem}\label[maintheorem]{thm:main}
Let $X$ be a positive-dimensional compact hyperbolic complex manifold.
Then $K_X$ is ample and $X$ is projective.
\end{maintheorem}

Let $X$ be a compact $C$-hyperbolic manifold. Then $X$ is hyperbolic since an entire curve lifts to the covering on which bounded
holomorphic functions separate points. Thus \Cref{thm:main} gives an affirmative answer to Kobayashi's question in \cite[Problem~7, p.~93]{Kob73}.
For a compact Riemann surface, \Cref{thm:main} follows from
\eqref{eq:curve-hyperbolicity}.  Wong proved that every compact
hyperbolic K\"ahler surface has ample canonical bundle
\cite[Theorem~3.1, p.~411]{Wong81}. For a smooth hyperbolic projective
surface $S$, the surface classification recalled by Green and Griffiths
\cite[\S4(a), pp.~67--68]{GG80} reduces the general-type assertion to
excluding projective K3 surfaces. These contain rational curves by the
Bogomolov--Mumford theorem (see Mori and Mukai \cite[Appendix]{MM83}).
Thus $S$ is of general type, and the absence of rational curves implies
that $K_S$ is ample \cite[Lemma~2.1]{DT}. In dimension two,
\Cref{thm:main} is also stated in \cite[Corollary~1]{QT26}
without the K\"ahler assumption.

In higher dimensions, Kodaira proved the projectivity for  compact quotients
of bounded domains by free, properly discontinuous groups of
biholomorphisms \cite[\S4]{Kod54}. OpenAI proved \Cref{thm:main} for all
compact hyperbolic K\"ahler manifolds \cite[Theorem~1.1]{OpenAI051}.
This result can be extended to all compact hyperbolic 
manifolds in Fujiki's class $\mathcal C$
\cite[Lemma~2.5]{QT26}, i.e., compact manifolds bimeromorphic to compact
K\"ahler manifolds.
Qin and Tian state that $K_X$ is ample if a compact hyperbolic manifold
$X$ carries a strictly positive real $(1,1)$-current $T$ such that $dT$ is smooth,
$dd^cT\le0$, and $\dim\{x\in X:\nu_x(T)\ge c\}\le2$ for every $c>0$.
Here $\nu_x(T)$ denotes the Lelong number of the local potential of $T$
after subtraction of a smooth $(1,1)$-form \cite[\S2 and Theorem~A]{QT26}. This includes the case of a
smooth Hermitian form $\omega$ with $dd^c\omega\le0$.
In particular, the SKT-hyperbolic case follows from the defining pluriclosed
metric and Marouani's hyperbolicity theorem \cite[Theorem~3.5]{Marouani}.
For divisorially hyperbolic surfaces, \cite[Proposition~3.1]{MarSurf26}
reduces canonical ampleness to the surface case above.
The smooth volume estimate in the proof of \Cref{thm:main} below
also gives the bigness step in \cite[\S3.4.1]{QT26}.

Let $h$ be a Hermitian metric on a compact complex manifold $X$, and let $R^h$
denote its Chern curvature tensor. The Chern holomorphic sectional
curvature (HSC) of $h$ is given by
\[
 H_h(\xi):=\frac{h(R^h(\xi,\bar\xi)\xi,\xi)}{|\xi|_h^4},
 \qquad x\in X,\quad 0\ne\xi\in T^{1,0}_xX.
\]
If $H_h<0$, continuity gives a constant $c>0$ such that
$H_h(\xi)\leq -c$ at every $x\in X$ and in every nonzero
direction $\xi\in T^{1,0}_xX$. The Schwarz lemma then implies that
$X$ is hyperbolic \cite[Theorem~3.8]{Kob67}.  Thus, \Cref{thm:main} gives the following corollary.

\begin{corollary}\label[corollary]{cor:negative-hsc}
Let $X$ be a compact complex manifold admitting a Hermitian metric $h$
with strictly negative Chern HSC. Then $K_X$ is ample and $X$ is
projective.
\end{corollary}

This corollary settles part~(c) of Conjecture~1.1 in \cite{YZ19}.
Suppose, in addition, that the Hermitian metric $h$ in
\Cref{cor:negative-hsc} is K\"ahler. Under this assumption, Heier,
Lu and Bun Wong gave a direct proof of ampleness for projective surfaces
and proved the result for projective threefolds
\cite[\S5.2 and Theorem~1.1]{HLW10}. Wu and Yau extended the result to projective
manifolds of arbitrary dimension \cite[Theorem~2]{WY}, and Tosatti
and Yang removed the projectivity assumption \cite[Corollary~1.2]{TY}.
Their respective preprints were first submitted to arXiv in May and June 2015.
Wu and Yau also gave a direct proof in the compact K\"ahler case
\cite[Theorem~1]{WYremark}. For K\"ahler metrics with quasi-negative HSC,
Pit-Mann Wong, Wu and Yau proved ampleness for projective manifolds of
Picard number one \cite{WWY12}. Diverio and Trapani
\cite[Theorem~1.2]{DT} and Wu and Yau \cite[Corollary~6]{WYremark}
proved ampleness in the general compact K\"ahler case. Their respective
preprints were first submitted to arXiv in June and September 2016.
Here quasi-negative HSC means that
\begin{equation}
  H_h\le 0, \qquad
  \exists\,x_0\in X \quad H_h(\xi)<0
      \quad\text{for every }0\ne\xi\in T^{1,0}_{x_0}X.
\end{equation}

On a compact K\"ahler manifold $X$, Broder and Stanfield proved that
$K_X$ is ample if $X$ admits a pluriclosed Hermitian metric $h$ ($\partial\bar\partial h=0$) with
strictly negative Chern HSC \cite[Theorem~1.1]{BS26}. Wan proved ampleness under quasi-negative Chern HSC
for compact K\"ahler manifolds and, more generally, for compact manifolds
in Fujiki's class $\mathcal C$ \cite[Theorem~1.1 and Remark~4.6]{Wan26}.

Suppose now that the metric $h$ in \Cref{cor:negative-hsc} has
constant Chern HSC $c<0$, meaning that $H_h(\xi)=c$ for every
$x\in X$ and every nonzero $\xi\in T^{1,0}_xX$. Since $X$ is
hyperbolic, \Cref{thm:main} implies that $X$ is projective and
hence belongs to Fujiki's class $\mathcal C$. If $X$ has complex
dimension at least two, then $h$ is K\"ahler by
\cite[Theorem~B(1)]{QT26}. 
\begin{corollary}
\label[corollary]{cor:intro-constant-hsc}
Let $h$ be a Hermitian metric of constant Chern HSC $c<0$ on a compact complex manifold $X$.
Then $h$ is K\"ahler and $(X,h)$ is a complex space form.
\end{corollary}
In complex dimension one, \Cref{cor:intro-constant-hsc} is immediate
since every Hermitian metric on a Riemann surface is K\"ahler.
The metric $h$ has constant HSC and is complete by compactness, so
$(X,h)$ is a complex space form. For completeness, we also give a direct proof of
\Cref{cor:intro-constant-hsc} in \S\ref{sec:ampleness}.

In 1977,  Wong stated without proof that a Hermitian metric of constant
nonzero HSC is K\"ahler \cite[p.~60]{Wong77}. By the classification of complex space forms
\cite[Theorem~I]{Hawley53}, \cite[Ch.~IX, Theorem~7.9]{KN69}, \Cref{cor:intro-constant-hsc}
answers his question whether the universal cover of a compact Hermitian
manifold with constant negative HSC is biholomorphic to the open unit ball
\cite[p.~353]{Wong91}.
Wong proved that a bounded domain $D\subset\mathbb C^n$ admitting such a compact
quotient by a group of biholomorphisms acting freely and properly discontinuously is biholomorphic to the open unit ball if there is a point $p\in\partial D$
through which no positive-dimensional complex analytic variety contained in
$\partial D$ passes and for which $U\cap D$ is biholomorphic to a bounded
convex domain for some neighbourhood $U$ of $p$ \cite[pp.~353--354]{Wong91}.

Balas and Gauduchon proved the K\"ahlerness conclusion for surfaces \cite[Theorem~1]{BG85}.
Qin and Tian proved the result for threefolds and for fourfolds with
balanced $h$ ($d(h^3)=0$ in dimension four) \cite[Theorems~C and~D(1)]{QT26}. Other results apply
under additional hypotheses on the metric. Chen, Chen and Nie treated
locally conformally K\"ahler metrics \cite[Theorem~1.1]{CCN21}, Rao and
Zheng treated Strominger K\"ahler-like metrics \cite{RZ22}, and Tang
treated pluriclosed metrics on compact K\"ahler manifolds
\cite[Theorem~1.5]{Tang26H}.

\subsection{Subvarieties of general type}\label{sec:comparison}
Let $M$ be a compact complex manifold of dimension $n$, with canonical bundle
$K_M:=\bigwedge^n(T^{1,0}M)^*$.  The canonical volume of $M$ is
\[
 \vol(K_M):=\limsup_{m\to\infty}\frac{n!}{m^n}\dim H^0(M,K_M^{\otimes m}).
\]
The manifold $M$ is of general type if and only if  $\vol(K_M)>0$.
For a reduced irreducible compact complex space $Z$, let $\widetilde Z$
be a resolution of its normalization. The space $Z$ is of general type
if $\vol(K_{\widetilde Z})>0$ \cite[\S1, p.~7]{KO75}. For a non-reduced
irreducible space $Z$ we apply this definition to its reduction $Z_{\rm red}$.
This definition is independent of the choice
of resolution, since any two such resolutions are bimeromorphic and therefore have the same plurigenera (the dimensions of the spaces of pluricanonical sections) and hence the same canonical volume
(\S\ref{sec:plurigenera}).

Recall that a reduced irreducible compact complex space $Z$ belongs to Fujiki's
class $\mathcal C$ if there is a surjective holomorphic map $M\to Z$
from a compact K\"ahler manifold $M$ \cite[\S4.3, Lemma~4.6(1)]{Fujiki78}. 

Let $g\colon Z\to X$ be a holomorphic map from a reduced irreducible
compact complex space $Z$. By the \emph{reduced image} of $g$ we mean
$Y=g(Z)$ equipped with the structure sheaf $(\mathcal O_X/\mathcal I_Y)|_Y$.
Here $\mathcal O_X$ is the structure sheaf of $X$, and $\mathcal I_Y$
is the ideal sheaf of local holomorphic functions on $X$ that vanish on $Y$.
Recall that $g$ is \emph{generically finite} if there
is a proper closed analytic subset $B\subsetneq Y$ such that
$g^{-1}(y)$ is finite for every $y\in Y\setminus B$.

The following result treats generically finite maps to compact hyperbolic complex spaces.

\begin{maintheorem}\label[maintheorem]{thm:proj}
Let $X$ be a compact hyperbolic complex space and let $Z$ be a
positive-dimensional reduced irreducible compact complex space in
Fujiki's class $\mathcal C$. If a holomorphic map $g\colon Z\to X$ is
generically finite, then $Z$ and the reduced image
$Y=g(Z)$ of $g$ are both of general type.
\end{maintheorem}

Together, \Cref{thm:main,thm:proj} give the following consequence.

\begin{corollary}\label[corollary]{cor:analytic-subspaces}
Let $X$ be a compact hyperbolic complex manifold. Every positive-dimensional
irreducible closed analytic subspace $Z\subseteq X$ is of general type.
\end{corollary}

Taking $g=\mathrm{id}_X$ in \Cref{thm:proj} proves the general-type
prediction of Kobayashi and Ochiai for positive-dimensional reduced
irreducible compact hyperbolic spaces in Fujiki's class $\mathcal C$
\cite[\S1, p.~7]{KO75}. Taking $g$ to be the inclusion of a
positive-dimensional irreducible closed subvariety of a hyperbolic
projective variety gives the forward implication of Lang's conjecture,
even when the ambient variety is singular
\cite[Conjecture~5.6]{Lang86}.

Hu, Meng and Zhang proved the conclusion of \Cref{thm:proj} when
$X$ and $Z$ are projective and $\dim Z\le 2$
\cite[Theorem~3.1]{HMZ}. In arbitrary dimension the projective image $Y$
is also of general type if its Albanese map is generically finite or its
Kodaira dimension is at least $\dim Y-2$ \cite[Proposition~1.7(i),(ii)]{HMZ}.
Suppose the reduced image $Y$ embeds in a
complex torus. Since $Y$ inherits hyperbolicity from $X$, Ueno's
structure theorem gives general type for $Y$
\cite[Theorem~3.10]{Ueno73}. More generally, if $Y$ admits a finite
holomorphic map to a complex torus, Kawamata's structure theorem
applied to the normalization $Y^\nu$ gives general type for $Y$
\cite[Theorem~23]{Kaw81Tor}. In each of these cases, general type for $Y$
implies general type for $Z$ since $g\colon Z\to Y$ is generically finite
\cite[Proposition~1.9]{Ueno73}.

For a positive-dimensional irreducible closed analytic subset $Y$ of a
compact K\"ahler manifold $X$, Guenancia proved that $Y$ is of general
type when $X$ carries a K\"ahler metric with strictly negative HSC
\cite[Theorem~A]{Gue}. Boucksom and Diverio obtained general type for
$Y$ when an \'etale cover of $X$ admits a bounded strictly
plurisubharmonic function \cite[Theorem~A(iv)]{BD21}. When the infinitesimal Carath\'eodory metric on the universal cover of $X$ is
nondegenerate, Kwok-Kin Wong and Yeung construct such a function as a
weighted sum of squared moduli of bounded holomorphic functions
\cite{WongYeung}.
Bei, Diverio, Eyssidieux and Trapani proved the inclusion case when
$X$ is K\"ahler hyperbolic in Gromov's sense, meaning that some
K\"ahler form on $X$ lifts to a form with a bounded primitive on the
universal cover of $X$, with boundedness measured in the lifted metric
\cite[Theorem~B]{BDET}. They also showed that every positive-dimensional
compact K\"ahler manifold admitting a generically finite holomorphic map
to such an $X$ is of general type \cite[Proposition~2.29(i) and Theorem~4.1]{BDET}.
Brunebarbe proved the inclusion case for hyperbolic
projective varieties carrying a large complex local system
\cite[Corollary~4.2]{Bru}. Here largeness means that the pullback to the
normalization of every positive-dimensional irreducible closed subvariety
has infinite monodromy \cite[Definition~3.1]{Bru}.

\subsection{Applications}

Golota proved canonical ampleness for projective manifolds with nondegenerate
negative total $k$-jet curvature \cite[Theorem~1.0.2]{Gol17}. For compact
K\"ahler manifolds, nondegenerate negative $k$-jet curvature already suffices
by Demailly's hyperbolicity theorem and the K\"ahler case of \Cref{thm:main}
\cite[Theorem~7.8]{Dem97}, \cite[Theorem~1.1]{OpenAI051}.
We obtain ampleness without a K\"ahler assumption, both for nondegenerate
negative $k$-jet curvature and for complex Finsler metrics of strictly negative
holomorphic sectional curvature (\Cref{cor:negative-directed-jets,cor:negative-finsler}).

For a projective target in \Cref{thm:proj}, the generic degree of the map
and the projective degree of its reduced image give a lower bound for the
canonical volume of every resolution of the source by \eqref{eq:volume-degree}.
The canonical bundle of the resolution need not be nef, an issue raised
in \cite[Question~1.3 and the discussion following it]{Zhang26}.
For curves, the genus--degree estimate already follows from the disc
comparison and Gauss--Bonnet \cite[proof of Theorem~2.1]{Dem97},
\cite[Theorem~1.2]{JK18}. In higher dimensions, when the source is also a projective manifold,
we obtain mixed intersection bounds on its relative canonical model (\Cref{cor:model-ample}).
We also obtain such bounds on the normalizations of the strata of a hyperbolic
projective semi-log canonical pair (\Cref{cor:slc-ample}).

Let $M\in\mathcal C$ be a compact complex manifold and $W$ a compact hyperbolic complex space.
Every meromorphic map $M\dashrightarrow W$ factors through Campana's core,
since the reduced image of any nonconstant such map is of general type.
If $M$ is weakly special, every meromorphic map from a connected finite
\'etale cover of $M$ to $W$ is constant. If $W$ is projective, every rational
map from an irreducible projective variety $S$ to $W$ factors through
Abramovich's Lang map of $S$
(\Cref{cor:weak-special-maps,cor:special-maps}).

\subsection{Proof strategy}\label{sec:proof-strategy}

\theoremstyle{plain}
\newcommand{\restatementheading}{Theorem}
\newtheorem*{restatementinner}{\restatementheading}
\newenvironment{restatement}[1]
 {\renewcommand{\restatementheading}{#1}\begin{restatementinner}}
 {\end{restatementinner}}
\theoremstyle{bolddefinition}

\introproofheading{\Cref{thm:main}.}

Let $h$ be a Hermitian metric on a compact hyperbolic complex manifold  $X$ of dimension $n$,
and let $\mathbb{D}$ be the open unit disc in $\mathbb{C}$.
A holomorphic disc in $X$ is a holomorphic map $f\colon\Disc\to X$.
Its center is $f(0)$ and its $h$-area is $A_h(f):=\int_\Disc f^*h$.
In \cite[\S\S3--5]{OpenAI051} and \cite[Propositions~3.6--3.7]{QT26},
one fixes a center $x_0\in X$ and maximizes a functional containing
the term $-\delta A_h(f)$, with $\delta>0$, over all finite-area
holomorphic discs centered at $x_0$. A disc attaining this supremum
is called a maximizing disc.
 At a holomorphic disc which attains this supremum, the argument requires an upper bound
for the second variation of $A_h$. This variation contains an
integral over $\Disc$ involving $dd^c h$. The integral
vanishes when $h$ is K\"ahler and is nonpositive when $\ddc h\le0$.
For a general Hermitian metric $h$, however, this integral has no prescribed sign.
We therefore replace the area term $-\delta A_h(f)$ by $-\delta\mathcal E(f)$, where $\mathcal E(f)$ is the Dirichlet energy of the
harmonic extension of the boundary trace of $I\circ f$ for a  smooth
embedding $I\colon X\hookrightarrow\R^N$. The area estimate $A_h(f)\le C\mathcal E(f)$ from \Cref{g:th:main}
and the second variation bound from \Cref{v:l:hess} allow us to derive
the pointwise volume comparison \eqref{eq:fixed-disc-volume} by varying
a maximizing disc. Taking $\beta=ah$ for large $a$ in this comparison
and using the constant disc give a bounded disc envelope. By
\cite[Theorem~1.2]{LSWZ26}, it yields a closed positive current in
$c_1(K_X)$ with bounded local potentials, hence zero Lelong numbers.
Thus $K_X$ is nef by \cite[Corollary~6.4]{DemReg92}.

Applying the pointwise volume comparison \eqref{eq:fixed-disc-volume} to
constant discs along the Hermitian continuity equation
\cite[Theorem~1.1]{ShermanWeinkove20} gives
$\beta_\varepsilon^n\ge((n+1)/D_h)^n h^n$, where
\[
D_h:=\sup_{f\colon\Disc\to X\ {\rm holomorphic}}|f'(0)|_h^2<\infty.
\] The mixed integral estimate \cite{CLN69} controls the
error in the strong holomorphic Morse inequality
\cite[Theorem~0.1(b)]{DemMorse85}, yielding
\begin{equation}\label{eq:strategy-canonical-volume}
 \vol(K_X)\ge\left(\frac{n+1}{2\pi D_h}\right)^n\int_X h^n>0 .
\end{equation}
Thus $K_X$ is big and $X$ is a Moishezon manifold \cite[Theorem~5.1]{DemMorse85}. Since $X$ contains no rational
curves, it is projective by \cite[Corollary~1.4.6]{BCHM10} (see \cite[Theorem~4.4]{HM10} for the  smooth Moishezon case). Since $X$ is of general type, projective and contains no ration curves, 
$K_X$ is ample by \cite[Lemma~2.1]{DT}.

\introproofheading{\Cref{thm:proj}.}
Let $Y$ be the reduced image of $g\colon Z\to X$ and put $n:=\dim Y=\dim Z$.
A resolution of $Z$ or $Y$ may contain rational curves.  Since $Y$ is hyperbolic,
we can choose a nonzero semipositive form $\tau$ supported in a coordinate
open ball in $Y_{\rm reg}$, positive on a smaller open ball, such that
$|q'(0)|_\tau^2\le1$ for every holomorphic disc $q$ in $Y$.
For a resolution $P$ of $Z$ or $Y$, let $F\colon P\to Y$ be the induced
map and put $\eta:=F^*\tau$. Then $|f'(0)|_\eta^2\le1$ for all holomorphic discs in $P$
and 
\begin{equation}
\int_P\eta^n=(\deg F)\int_{Y_{\rm reg}}\tau^n>0.
\end{equation}

Take a K\"ahler modification $\rho\colon Q\to P$ and a K\"ahler metric $h$
on $Q$. The derivative bound passes to $\rho^*\eta$, which may degenerate.
For holomorphic discs in $Q$, the sum of the boundary energy $\mathcal E(f)$ and the
Green energy $Q_h(f)$ controls $A_h(f)$, and bounded Green energy prevents
interior area concentration along a sequence of bounded area. The disc variation now gives
$\beta^n\ge(n+1)^n(\rho^*\eta)^n$ whenever $\tr_\beta\Ric^{(1)}_\beta\ge-n$.
Testing this inequality against Gauduchon metrics gives a degree bound for
curves. In higher dimensions, it gives a positive current $T\in c_1(K_Q)$
with $(2\pi)^n T_{\rm ac}^n\ge(n+1)^n(\rho^*\eta)^n$ almost everywhere,
where $T_{\rm ac}$ is the absolutely continuous part of $T$. Since $Q$ is K\"ahler, Boucksom's estimate
gives $\vol(K_Q)\ge\int_Q T_{\rm ac}^n$ \cite[Proposition~3.1]{Bou02}.
Since $\deg\rho=1$ and canonical volume is bimeromorphically invariant,
in both cases we obtain
\begin{equation}\label{eq:strategy-relative-volume}
 \vol(K_P)\ge\left(\frac{n+1}{2\pi}\right)^n\int_P\eta^n>0 .
\end{equation}
Taking $F$ first from a resolution of $Z$ and then from a resolution of $Y$
gives general type for $Z$ and $Y$, respectively.

\section{Preliminaries}\label{sec:preliminaries}
In this section, we firstly fix some notation and conventions that will be used throughout. 
We then establish the integral identities, convergence of boundary traces,
and uniform derivative bound needed in \S\ref{sec:disc-estimates}.

\subsection{Notation and conventions}\label{sec:conventions}

Let $X$ be a topological space, and let $A\subset X$ be a subset. Then we use $\overline{A}$ to denote the closure of 
$A$. For any subset $U\subset X$,  we write $A\Subset U$ if $\overline{A}\subset U$ is a compact subset.  For any $x\in \mathbb{R}^d$, and 
any $r>0$, set $B(x,r):=\{y\in\mathbb{R}^d: |y-x|<r\}.$ 
Write $\Disc=\{z\in\C:|z|<1\}$, $\overline\Disc=\{z\in\C:|z|\le1\}$
and $\T=\partial\Disc=\{z\in\C:|z|=1\}$.
For $z=a+\sqrt{-1}b$, put
$\partial_z=(\partial_a-\sqrt{-1}\partial_b)/2$ and
$\partial_{\bar z}=(\partial_a+\sqrt{-1}\partial_b)/2$.
We define the Euclidean Laplacian on $\mathbb{C}$ by $\Delta_{\rm eucl}:=\partial_a^2+\partial_b^2$.
Let $\lambda$ be the planar Lebesgue measure, so $d\lambda=da\,db,$ 
$\lambda(\Disc)=\pi$. We always identify a real two-form on a planar domain with the signed
measure it induces in the standard orientation. In particular,
$(\sqrt{-1}/2)\,dz\wedge d\bar z=da\wedge db$ induces $\lambda$.
On $\T$, write $\zeta=e^{\sqrt{-1}\theta}$ and use
$d\sigma=d\theta/(2\pi)$ and $\av_\T u=\int_\T u\,d\sigma$.

In $C^\infty(\overline U)$ all derivatives extend continuously to
$\overline U$. For scalar functions
or maps into a finite-dimensional real or complex inner-product space $V$,
\[
 \|u\|_{L^p(S,\mu)}=\left(\int_S|u|^p\,d\mu\right)^{1/p}\quad(1\le p<\infty),
 \qquad \|u\|_{L^\infty(S,\mu)}=\operatorname*{ess\,sup}_S |u| .
\]
We abbreviate these norms to $\|u\|_p$ when the domain and measure are
clear, using $d\lambda$ on $\Disc$ and $d\sigma$ on $\T$.
Using componentwise weak derivatives, we take
\[
 \begin{aligned}
 \|u\|_{W^{1,p}(S;V)}^p&=\|u\|_p^p+\|\partial_a u\|_p^p+\|\partial_b u\|_p^p\quad(p=2,4),\\
 \|u\|_{W^{1,2}(\T;V)}^2&=\|u\|_2^2+\|\partial_\theta u\|_2^2 ,
 \end{aligned}
\]
where $S\subset\C$ is a domain and $\partial_\theta$ is the periodic
angular derivative. The zero-trace space $W^{1,2}_0(\Disc;V)$ is the
closure of $C_c^\infty(\Disc;V)$ in $W^{1,2}(\Disc;V)$.
For $S=\Disc$ or $\T$ and a smooth manifold $X$, we test membership
and convergence in $W^{1,p}(S,X)$ after composing with a smooth embedding
$I\colon X\hookrightarrow\R^N$ by Whitney's embedding theorem \cite[\S8, Theorem~1]{Whitney36}. 

For $F\colon\Disc\to\C^N$ we use the real Hilbert--Schmidt norm
\begin{equation}\label{eq:disc-hilbert-schmidt}
 |dF|^2=|\partial_aF|^2+|\partial_bF|^2
       =\sum_{j=1}^N\bigl(|\partial_a F^j|^2+|\partial_b F^j|^2\bigr).
\end{equation}
For Sobolev maps this identity holds almost everywhere.
For $u\in L^2(\T;V)$ and $m\in\mathbb Z$, put
$\widehat u(m):=\int_\T u(\zeta)\bar\zeta^m\,d\sigma(\zeta)$,
using the complexification of $V$ with its induced Hermitian norm when $V$ is real.
The boundary space $H^{1/2}(\T;V)$ has squared norm
\begin{equation}\label{h1/2}
 \|u\|_{H^{1/2}(\T;V)}^2=\sum_{m\in\mathbb Z}(1+|m|)\,|\widehat u(m)|^2 .
\end{equation}
For holomorphic $u$ on $\Disc$ the Hardy norm is
$\|u\|_{H^2(\Disc)}^2=\sup_{0<r<1}\int_\T|u(r\zeta)|^2\,d\sigma(\zeta)$.

All complex manifolds are connected, smooth and of positive dimension unless
otherwise stated. A projective complex space is one admitting a closed holomorphic embedding into some $\mathbb P^N$.
For rational maps from a complex algebraic variety to a projective complex space $W$,
we identify $W_{\rm red}$ with its algebraic model by Chow's theorem \cite[Theorem~V]{Chow49}.
On smooth forms of every degree we use the normalization of \cite[\S9.B]{Dem20},
\[
 d^c:=\frac{\sqrt{-1}}{4\pi}(\bar\partial-\partial),\qquad
 \ddc=d\circ d^c=\frac{\sqrt{-1}}{2\pi}\partial\bar\partial .
\]
These operators act distributionally on currents. 
Let $X$ be a compact complex manifold of dimension $n$,  write $T_X=T^{1,0}X$ for the holomorphic tangent
bundle and $K_X=\bigwedge^n T_X^*$ for the canonical bundle, so the corresponding first Chern classes satisfy
$c_1(X)=c_1(T_X)=-c_1(K_X)$. A compact complex manifold $X$
is of general type when its canonical bundle $K_X$ is big, equivalently
when $\vol(K_X)>0$ (\S\ref{sec:introduction}).
All Hermitian metrics are smooth unless otherwise stated.
For a Hermitian metric $h$ on $T_X$, we take $h(v,w)$ to be linear in $v$ and identify $h$ with its fundamental $(1,1)$-form.
For $v\in T^{1,0}_xX$, put $|v|_h:=\sqrt{h_x(v,v)}$. The associated
real Riemannian metric is $h_{\R}(u,v):=2\operatorname{Re}h(u^{1,0},v^{1,0})$.
We use $\dis_h$ and $\operatorname{Length}_h$ for its distance and length.
For a Hermitian metric $h$ on a compact complex manifold, its volume form is given by $h^n/n!$.

For a positive Hermitian form $h$ and a real $(1,1)$-form $\alpha$,
in local holomorphic coordinates, write 
\[
 h=\sqrt{-1}\sum_{j,k}h_{j\bar k}\,dz^j\wedge d\bar z^k,\qquad
 \alpha=\sqrt{-1}\sum_{j,k}a_{j\bar k}\,dz^j\wedge d\bar z^k.
\]
The trace $\tr_h\alpha$ is the contraction of $\alpha$ with the inverse
metric of $h$. With $H=(h_{j\bar i})_{i,j}$ and $A=(a_{j\bar i})_{i,j}$, we have
\begin{equation}\label{eq:hermitian-trace-wedge}
 \tr_h\alpha=\tr(H^{-1}A),\qquad
 \alpha\wedge h^{n-1}=\frac{\tr_h\alpha}{n}\,h^n,
\end{equation}
where $\tr(\cdot)$ denotes the matrix trace. For tangent vectors written as columns, $h(v,w)=w^*Hv$, with
$w^*$ the conjugate transpose. We put $\alpha[v]:=-\sqrt{-1}\,\alpha(v,\bar v)=v^*Av$
and $|v|_\alpha^2:=\alpha[v]$ when $\alpha\ge0$. We write $\alpha\ge0$ (respectively $\alpha>0$) when $A$ is
positive semidefinite (respectively positive definite), and order real $(1,1)$-forms by their difference.
On smooth functions $u$, put $\Delta_h u:=2\pi\tr_h(\ddc u)$.
The first Chern--Ricci form of $h$ is $\Ric^{(1)}_h:=-2\pi\ddc\log\det H$,
and its trace $\tr_h\Ric^{(1)}_h$ is the Chern scalar curvature of $h$.
The Chern connection $\nabla^h$ on $T_X$ preserves $h$ and has $(0,1)$
part $\bar\partial$. Its curvature tensor and the Chern holomorphic
sectional curvature are
\[
 R^h(U,V):=[\nabla^h_U,\nabla^h_V]-\nabla^h_{[U,V]},\qquad
 H_h(\xi):=\frac{h(R^h(\xi,\bar\xi)\xi,\xi)}{|\xi|_h^4}
 \quad(0\ne\xi\in T^{1,0}_xX),
\]
where $[\cdot,\cdot]$ is Lie bracket.

Let $h$ be a Hermitian metric on a compact complex manifold $X$ of dimension $n$.
The real Bott--Chern group of bidegree $(1,1)$ consists of smooth
$d$-closed real $(1,1)$-forms modulo $\ddc u$ for smooth real functions $u$ on $X$.
A class $\alpha$ in this group is nef if, for every $\varepsilon>0$, it admits a
smooth representative $\alpha_\varepsilon$ satisfying $\alpha_\varepsilon\ge-\varepsilon h$.
This condition is independent of $h$ since any two Hermitian metrics on compact $X$
bound each other by constant multiples.

\subsection{Holomorphic discs}\label{sec:disc-conventions}\label{sec:boundary-traces}

In this subsection, we collect some estimates on the unit disc for later use.
We begin with a circular mean identity on the unit disc.
\begin{lemma}\label[lemma]{lem:disc-jensen-flux}
Let $u\in C^2(\Disc)\cap C^0(\overline\Disc)$ be a real-valued function such that 
$\Delta_{\rm eucl}u\in L^1(\Disc)$. Then
\begin{equation}\label{eq:disc-jensen}
 \int_\Disc \log|z|^{-2}\,\ddc u=\av_\T u-u(0).
\end{equation}
If, in addition, $u\ge0$ on $\Disc$ and $u|_\T=0$, then
\begin{equation}\label{eq:boundary-flux-sign}
 \int_\Disc\Delta_{\rm eucl}u\,d\lambda\le0 .
\end{equation}
\end{lemma}
\begin{proof}
For $0<r<1$, set $m(r):=\av_\T u(r\,\cdot)$. 
Since $u\in C^2(\Disc)$, the Green--Riesz representation formula \cite[pp.~356--357]{Riesz30}, in the normalization of \cite[Ch.~I, \S4.B.1, (4.5)--(4.6)]{DemCADG},
gives
\begin{equation}\label{eq:disc-jensen-radius}
 m(r)-u(0)=\int_{r\Disc}\log\frac{r^2}{|z|^2}\,\ddc u
 \qquad(0<r<1).
\end{equation}
Letting $r\uparrow1$, the dominated convergence theorem and continuity of $u$
give \eqref{eq:disc-jensen}.

By the divergence theorem,
\begin{equation}\label{eq:disc-flux-radius}
 2\pi r m'(r)=\int_{r\Disc}\Delta_{\rm eucl}u\,d\lambda.
\end{equation}
If $\int_\Disc\Delta_{\rm eucl}u\,d\lambda>0$, then $m$ is strictly
increasing near $1$, contrary to $m(r)\ge0=m(1)$.
\end{proof}

A \emph{holomorphic disc} in a complex space $X$ is a holomorphic map
$g\colon\Disc\to X$. The point $g(0)\in X$ is called the \emph{center} of $g$.
A closed holomorphic disc in $X$
is a holomorphic disc $g\colon\Disc\to X$ that can be extended holomorphically
to an open neighbourhood of $\overline\Disc$. Now let $(X,h)$ be a compact
Hermitian manifold. Writing $g'(z):=dg_z(\partial_z)\in T^{1,0}_{g(z)}X$, we have 
\begin{equation}\label{eq:disc-local-identities}
 g^*h=2|g'(z)|_h^2\,d\lambda .
\end{equation}
By Whitney's embedding theorem \cite[\S8, Theorem~1]{Whitney36}, there is a smooth embedding
$I\colon X\hookrightarrow\R^N$. In the following, we always fix such an embedding $I$.

For any holomorphic disc $g\colon\Disc\to X$, define 
\begin{equation}
 A_h(g):=\int_\Disc g^*h
       =2\int_\Disc|g'(z)|_h^2\,d\lambda(z)\in[0,+\infty].
\end{equation}
Then we call $g$ a \emph{finite-area holomorphic disc} if $A_h(g)<\infty$.
This definition is independent of the choice of $h$ as $X$ is compact.
For $F:=I\circ g$, put $a_I:=2\min_{|v|_{h_{\R}}=1}|dI(v)|^2$ and
$b_I:=2\max_{|v|_{h_{\R}}=1}|dI(v)|^2$. Compactness and injectivity of $dI$
give $0<a_I\le b_I$. Since $|dg(\partial_a)|_{h_{\R}}^2=|dg(\partial_b)|_{h_{\R}}^2=2|g'|_h^2$,
\begin{equation}\label{eq:embedded-disc-energy-density}
 2a_I|g'|_h^2\le |dF|^2\le 2b_I|g'|_h^2 .
\end{equation}
If $A_h(g)<\infty$, then $F\in W^{1,2}(\Disc;\R^N)$.
For a finite-dimensional real or complex inner-product space $V$, the trace theorem
\cite[pp.~289--290, (1.4)--(1.5)]{Gagliardo57}, \cite[Theorems~3.37 and~3.40(ii)]{McLean00} applied componentwise gives a unique bounded linear extension of restriction from $C^\infty(\overline\Disc,V)$,
\[
 \operatorname{Tr}_\T\colon W^{1,2}(\Disc;V)\longrightarrow H^{1/2}(\T;V),
 \qquad \ker\operatorname{Tr}_\T=W^{1,2}_0(\Disc;V).
\] If $g$ is a finite area holomorphic disc, then the
\emph{boundary trace} $\gamma_g$ of  $g$ is defined by
$\gamma_g:=\operatorname{Tr}_\T(I\circ g)\in H^{1/2}(\T;\R^N)$.
If $g$ extends continuously to $\overline\Disc$, then $\gamma_g=I\circ(g|_\T)$.

By the Poisson kernel of the unit disc\cite[Ch.~I, \S4.B, (4.6)--(4.7)]{DemCADG},
we define the Poisson extension of $u\in L^2(\T;V)$
and its Dirichlet energy by
\[
 (\mathsf Pu)(z):=\int_\T\frac{1-|z|^2}{|\zeta-z|^2}\,u(\zeta)\,d\sigma(\zeta),
 \qquad H(u):=\int_\Disc|d(\mathsf Pu)|^2\,d\lambda .
\]
With $d\sigma=d\theta/(2\pi)$, the Fourier representation of
$\mathsf Pu$ and Douglas's formula
\cite[\S5, (5.1), and \S21, (21.1)--(21.2)]{Douglas}\cite[Lemma~1.5.2]{EKMR14} give
\begin{equation}\label{g:eq:douglas}
 H(u)=2\pi\sum_{m\in\mathbb Z}|m|\,|\widehat u(m)|^2
 =2\pi\int_\T\int_\T
 \frac{|u(\zeta)-u(\eta)|^2}{|\zeta-\eta|^2}\,
 d\sigma(\zeta)\,d\sigma(\eta).
\end{equation}
The \emph{boundary energy} of a finite-area holomorphic disc $g$ is defined by 
$\mathcal E(g):=H(\gamma_g)$. By \eqref{h1/2} and \eqref{g:eq:douglas},
$H(u)<\infty$ for $u\in L^2(\T;V)$ exactly when $u\in H^{1/2}(\T;V)$, so $\mathcal E(g)<\infty$.

For $u\in H^{1/2}(\T;V)$, the Poisson extension $\mathsf Pu$ belongs
to $W^{1,2}(\Disc;V)$ and has boundary trace $u$
\cite[\S2.1, (2.3), and Lemma~A.5]{VW20}.
Dirichlet's principle
\cite[\S8.2.2, Theorem~3; \S8.2.3, (38) and the remark after Theorem~4,
pp.~471--475]{Evans10} gives
$\int_\Disc\langle d(\mathsf Pu),dv\rangle\,d\lambda=0$
for every $v\in W^{1,2}_0(\Disc;V)$. It also gives that
$\mathsf Pu$ uniquely minimizes the Dirichlet energy among maps in
$W^{1,2}(\Disc;V)$ with boundary trace $u$. 

Recall that a sequence of smooth maps $g_j\colon\Disc\to X$ converge locally smoothly to a smooth map $g\colon\Disc\to X$
if, for every $0<r<1$, $g_j\to g$ uniformly with all partial derivatives on $\{|z|\le r\}$ in local coordinates on $X$.

\begin{lemma}
\label[lemma]{lem:radial-traces}
Let $(X,h)$ be a compact Hermitian manifold with a smooth embedding
$I\colon X\hookrightarrow\R^N$. Then
\begin{enumerate}
\renewcommand{\theenumi}{\thetheorem(\roman{enumi})}
\renewcommand{\labelenumi}{\textup{(\roman{enumi})}}

\item\label[lemma]{lem:radial-approximation}
Let $g\colon \Disc\to X$ be a finite-area holomorphic disc and put $g_r(z):=g(rz)$
for $z\in\Disc$ and $0<r<1$. Then, as $r\uparrow1$,
\begin{equation}\label{eq:radial-disc-limits}
\begin{aligned}
 A_h(g_r)&\uparrow A_h(g),\qquad \gamma_{g_r}\longrightarrow\gamma_g
       \quad\text{in }H^{1/2}(\T;\R^N),\qquad  \mathcal E(g_r)\longrightarrow\mathcal E(g).
\end{aligned}
\end{equation}

\item\label[lemma]{lem:bounded-area-limits}
Let $g_j\colon\Disc\to X$ be holomorphic discs converging locally smoothly
to a map $f\colon \Disc\to X$, with $\sup_j A_h(g_j)<\infty$. Then $f$ is a holomorphic disc and 
\begin{equation}\label{eq:bounded-area-disc-limits}
\begin{aligned}
 A_h(f)&\le\liminf_{j\to\infty}A_h(g_j)<\infty,\qquad 
 \gamma_{g_j}\longrightarrow\gamma_f
       \quad\text{in }L^2(\T;\R^N),\qquad 
 \mathcal E(f)\le\liminf_{j\to\infty}\mathcal E(g_j).
\end{aligned}
\end{equation}
Moreover, for every smooth real $(1,1)$-form $\alpha$ on $X$,
\begin{equation}\label{eq:disc-weighted-pullback-limit}
 \lim_{j\to\infty}\int_\Disc \log|z|^{-2}\,g_j^*\alpha
 =\int_\Disc \log|z|^{-2}\,f^*\alpha .
\end{equation}
\end{enumerate}
\end{lemma}
\begin{proof}
\textup{(i)} It is clear that 
\begin{equation}
 A_h(g_r)=\int_{r\Disc}g^*h\uparrow A_h(g)
 \qquad\text{as }\qquad r\uparrow1.
\end{equation}
Set $F:=I\circ g$ and $F_r(z):=F(rz),\ z\in\Disc$. By
\eqref{eq:embedded-disc-energy-density}, $F\in W^{1,2}(\Disc;\R^N)$.
Since
\begin{equation}
 \|F_r\|_{L^2(\Disc)}\le r^{-1}\|F\|_{L^2(\Disc)},\qquad
 \|dF_r\|_{L^2(\Disc)}=\|dF\|_{L^2(r\Disc)},
\end{equation}
 $F_r$ converges to $F$ in $W^{1,2}(\Disc;\R^N)$ by density
\cite[\S5.3.3, Theorem~3]{Evans10}. The continuity of the trace operator $\operatorname{Tr}_\T$ then gives
\begin{equation}
 \gamma_{g_r}=F(r\,\cdot)|_\T\longrightarrow\gamma_g
 \qquad\text{in }H^{1/2}(\T;\R^N)\qquad\text{as }\qquad r\uparrow1.
\end{equation}
Thus, by \eqref{h1/2} and \eqref{g:eq:douglas},
\begin{equation}
\mathcal E(g_r)=H(\gamma_{g_r})\to H(\gamma_g)=\mathcal E(g).
\end{equation}
\smallskip
\textup{(ii)} By the assumption,  $f$ is a holomorphic disc, and for any 
$0<r<1$, 
\begin{equation}
 \int_{r\Disc}f^*h
 =\lim_{j\to\infty}\int_{r\Disc}g_j^*h
 \le\liminf_{j\to\infty}A_h(g_j).
\end{equation}
Set $A_*:=\sup_j A_h(g_j)<\infty$, then by letting $r\uparrow 1$,
\begin{equation}
A_h(f)\le\liminf_{j\to\infty}A_h(g_j)\le A_*.
\end{equation}

For any holomorphic disc
$q\colon\Disc\to X$, set $U:=I\circ q$. For $0<r<R<1$ and $\zeta\in\T$,
estimating the radial integral in $L^2(s\,ds)$ gives
\begin{equation}
 |U(R\zeta)-U(r\zeta)|^2
 =\left|\int_r^R\partial_s U(s\zeta)\,ds\right|^2
 \le \log(R/r)\int_r^R|\partial_s U(s\zeta)|^2s\,ds.
\end{equation}
Integration over any arc $J\subseteq\T$ then gives
\begin{equation}\label{eq:radial-trace}
 \int_J |U(R\zeta)-U(r\zeta)|^2\,d\sigma(\zeta)
 \le \frac{\log(R/r)}{2\pi}
     \int_{\{s\zeta:\ r<s<R,\ \zeta\in J\}} |dU|^2\,d\lambda .
\end{equation}

Write $F_j:=I\circ g_j$, $F:=I\circ f$.
Apply \eqref{eq:radial-trace} with $J=\T$ to $F_j$ and $F$, and let
$R\uparrow1$. By \textup{(i)}, \eqref{eq:embedded-disc-energy-density},
and the triangle inequality,
\begin{equation}\label{eq:radial-boundary-trace-comparison}
 \|\gamma_{g_j}-\gamma_f\|_{L^2(\T)}
 \le 2\sqrt{\frac{b_I A_*}{2\pi}\log(1/r)}
      +\|F_j(r\,\cdot)-F(r\,\cdot)\|_{L^2(\T)} .
\end{equation}
Letting first $j\to\infty$ and then $r\uparrow1$ gives
$\gamma_{g_j}\to\gamma_f$ in $L^2(\T;\R^N)$ as $j\to \infty$. Their Fourier coefficients
therefore converge at each index. Passing to the limit in finite partial sums of
\eqref{g:eq:douglas} and taking their supremum gives
\begin{equation}
\mathcal E(f)\le\liminf_{j\to\infty}\mathcal E(g_j).
\end{equation}

Finally, fix any   smooth real $(1,1)$-form $\alpha$ on $X$. On each
closed subdisc $\overline{r\Disc}$, the densities of $g_j^*\alpha$
converge uniformly to that of $f^*\alpha$. Since $\log|z|^{-2}$ is
integrable there, the weighted integrals over $r\Disc$ converge.
With $C_\alpha:=\max_{|v|_h=1}|\alpha[v]|$ we have
$-C_\alpha h\le\alpha\le C_\alpha h$. Thus, for $q=g_j$ or $q=f$,
\begin{equation}
 \left|\int_{\Disc\setminus r\Disc}\log|z|^{-2}\,q^*\alpha\right|
 \le C_\alpha A_*\log(r^{-2}).
\end{equation}
Note that this bound tends to zero as $r\uparrow1$, uniformly in $j$, then 
\eqref{eq:disc-weighted-pullback-limit} follows.
\end{proof}

\subsection{Derivative bounds of holomorphic discs}\label{sec:brody-derivatives}
The following derivative bound uses Brody's reparametrization lemma
\cite{Brody78}.
\begin{lemma}\label[lemma]{le:brody-derivative-bound}
Let $(X,h)$ be a compact hyperbolic Hermitian manifold. The supremum $D_h:=\sup_{f\colon\Disc\to X\ {\rm holomorphic}}|f'(0)|_h^2$
is positive and finite, and every
holomorphic disc $g\colon\Disc\to X$
satisfies
\begin{equation}\label{g:eq:brody-speed}
 (1-|z|^2)^2|g'(z)|_h^2\le D_h,\qquad z\in\Disc .
\end{equation}
\end{lemma}
\begin{proof}
Fix $x_0\in X$, choose a local holomorphic chart $(V,\varphi)$ around $x_0$ with $\varphi(x_0)=0$. For sufficiently
small $\varepsilon>0$, the holomorphic disc $g(z):=\varphi^{-1}(\varepsilon z,0,\ldots,0)$
satisfies $g'(0)\ne0$. Hence $D_h\ge |g'(0)|_h^2>0$. 

If $D_h=\infty$, Brody's lemma \cite{Brody78}, \cite[\S1, p.~1]{Duval21} gives a
nonconstant entire curve $\C\to X$, contradicting the hyperbolicity
of $X$. Thus $D_h<\infty$.

Composing $g$ with
$\phi_z(w)=(z+w)/(1+\bar z w)$ gives
$|(g\circ\phi_z)'(0)|_h=(1-|z|^2)|g'(z)|_h$. Then
\eqref{g:eq:brody-speed} follows by definition.
\end{proof}

\section{Estimates for holomorphic discs}\label{sec:disc-estimates}

We first bound the area of a holomorphic disc by its boundary energy.
We then vary a maximizing disc to compare Hermitian volume forms at its center.

\subsection{An area estimate for holomorphic discs}\label{sec:geometry}

\begin{theorem}[Area and boundary energy]\label[theorem]{g:th:main}
Let $(X,h)$ be a compact hyperbolic Hermitian manifold and fix a smooth
embedding $I\colon X\hookrightarrow\R^N$. There is a constant $C=C(X,h,I)\ge1$ such that
every finite-area holomorphic disc $g\colon\Disc\to X$ satisfies
\begin{equation}\label{g:eq:area-boundary}
 C^{-1}\mathcal E(g)\le A_h(g)\le C\mathcal E(g).
\end{equation}
\end{theorem}

Note that the hyperbolicity assumption in \Cref{g:th:main} is essential for the
upper bound in \eqref{g:eq:area-boundary}, which can fail even on $\mathbb P^1$.
In the affine coordinate $w$, take $h=\sqrt{-1}\,dw\wedge d\bar w/(1+|w|^2)^2$,
$I(w)=(2w,|w|^2-1)/(1+|w|^2)\in\C\times\R$, and $g_R(z)=Rz$.
The area integral and the boundary Fourier modes give
\[
 A_h(g_R)=\frac{2\pi R^2}{1+R^2},\qquad
 \mathcal E(g_R)=\frac{8\pi R^2}{(1+R^2)^2},\qquad
 \frac{A_h(g_R)}{\mathcal E(g_R)}=\frac{1+R^2}{4}\longrightarrow\infty
 \quad(R\to\infty).
\]

The negative Dirichlet Green function of $\Disc$ is
\cite[Ch.~I, \S4.A, (4.1)--(4.3)]{DemCADG}
\begin{equation}\label{eq:dirichlet-kernel-at-centre}
 \mathcal{G}_{\Disc}(z,w)=\frac1{2\pi}\log\left|\frac{z-w}{1-\bar wz}\right|\le0,
 \qquad \Delta_{{\rm eucl},z}\mathcal{G}_{\Disc}(\cdot,w)=\delta_w ,
\end{equation}
where $\delta_w$ is the Dirac measure concentrated at $w$. Thus $-4\pi\mathcal{G}_{\Disc}(z,0)=\log|z|^{-2}$.
For a holomorphic disc $g$ in $(X,h)$, define the \emph{Green energy} of
its area measure $\mu_g=g^*h$ by
\begin{equation}\label{eq:green-energy}
 Q_h(g):=-2\iint_{\Disc^2}\mathcal{G}_{\Disc}(z,w)\,d\mu_g(z)\,d\mu_g(w)\in[0,+\infty].
\end{equation}

\begin{lemma}\label[lemma]{g:le:green-hardy}
Let $(X,h)$ be a compact Hermitian manifold with a smooth embedding $I\colon X\to\mathbb{R}^N$, then there is a constant $C:=C(X,h,I)\ge1$ such that for any  finite-area holomorphic
disc $g$ in $X$,
\begin{equation}\label{g:eq:scalar-green}
 A_h(g)\le C\bigl(\mathcal E(g)+Q_h(g)\bigr),\qquad
 Q_h(g)\le8\bigl\|(1-|z|^2)|g'|_h^2\bigr\|_2^2 .
\end{equation}
\end{lemma}
\begin{proof}
Write $\rho_g(z):=|g'(z)|_h^2$ and $\chi(z):=1-|z|^2, z\in\Disc$, so that
$d\mu_g=2\rho_g\,d\lambda$ by \eqref{eq:disc-local-identities}.

For $0<r<1$, put $D_r:=r\Disc.$ 
Standard elliptic theory gives a unique real-valued solution
$p_r\in C^\infty(\overline{D_r})$ of the Dirichlet problem
\begin{equation}
 -\Delta_{\rm eucl} p_r=2\rho_g \quad\text{ in } D_r,\qquad
 p_r=0 \quad\text{ on } \partial D_r .
\end{equation}
Extend $p_r$ by zero from $D_r$ to $\Disc$. Its zero boundary trace gives $p_r\in W^{1,2}_0(\Disc)$, so the Hardy
inequality \cite[Theorem~3.2(1), (3.2)]{Avkhadiev25} gives
\begin{equation}\label{eq:weighted-hardy}
 \|p_r/\chi\|_2\le\|dp_r\|_2 .
\end{equation}
Integration by parts on $D_r$ and the $L^2$ pairing give
\begin{equation}
 \|dp_r\|_2^2=2\int_{D_r}\rho_g p_r\,d\lambda
 \le 2\|\chi \rho_g\|_{L^2(D_r)}\|p_r/\chi\|_2
 \le 2\|\chi \rho_g\|_{L^2(D_r)}\|dp_r\|_2 .
\end{equation}
Hence 
\begin{equation}\label{energy1}
\|dp_r\|_2^2\le4\|\chi \rho_g\|_{L^2(D_r)}^2.
\end{equation}

Let $g_r(z):=g(rz)$ and write $G_r(z,w):=\mathcal G_\Disc(z/r,w/r).$
Then the Green--Riesz formula \cite[pp.~356--357]{Riesz30}, \cite[Ch.~I, \S4.B.1, (4.5)]{DemCADG} gives
\begin{equation}
p_r(z)=-\int_{D_r}G_r(z,w)\,d\mu_g(w).
\end{equation}
 Changing variables in
\eqref{eq:green-energy} therefore gives
\begin{equation}\label{enery2}
Q_h(g_r)=2\int_{D_r}p_r\,d\mu_g=2\|dp_r\|_2^2.
\end{equation}
Extend $-G_r$ by zero outside $D_r^2$. By domain monotonicity of the Dirichlet Green function, $-G_r$ increases to $-\mathcal G_\Disc$ off the diagonal as $r\uparrow1$.
Hence the monotone convergence theorem gives
\begin{equation}\label{eq:contracted-green-kernel}
 Q_h(g_r)\uparrow Q_h(g)\qquad (r\uparrow1).
\end{equation}
From \eqref{energy1}, \eqref{enery2} and \eqref{eq:contracted-green-kernel}, we have
\begin{equation}\label{eq:green-energy-weighted-bound}
Q_h(g)\le8\|\chi \rho_g\|_2^2,
\end{equation}
which is  the second inequality in \eqref{g:eq:scalar-green}.

By part~\textup{(i)} of \Cref{lem:radial-traces},
\begin{equation}\label{eq:radial-area-boundary-limits}
 A_h(g_r)\uparrow A_h(g),\qquad
 \mathcal E(g_r)\to\mathcal E(g)\qquad (r\uparrow1).
\end{equation}
Together with \eqref{eq:contracted-green-kernel}, 
it suffices to prove the first inequality in \eqref{g:eq:scalar-green}
for closed holomorphic discs with a constant independent of the disc.

 Now we assume  $g$ is a closed holomorphic disc and put $F:=I\circ g$ and
$v:=F-\mathsf P\gamma_g$. Since $v|_\T=0$, $\mathsf P\gamma_g$ is
harmonic in $\Disc$, integration by parts gives
\begin{equation}\label{eq:disc-harmonic-orthogonality}
 \int_\Disc\langle d(\mathsf P\gamma_g),dv\rangle\,d\lambda=0.
\end{equation}
By \eqref{eq:embedded-disc-energy-density},
\begin{equation}\label{eq:disc-harmonic-energy-decomposition}
 a_I A_h(g)\le\|dF\|_2^2
 =\|d(\mathsf P\gamma_g)\|_2^2+\|dv\|_2^2
 =\mathcal E(g)+\|dv\|_2^2 .
\end{equation}
Let $\partial I$ be the $(1,0)$-part of the complexified differential of $I$,
and take $\partial\bar\partial I$ componentwise. Since $g$ is holomorphic,
$\Delta_{\rm eucl}v=4(\partial\bar\partial I)_{g(z)}(g',\bar g')$.
With $L_I:=\max\{1,4\max_{|v|_h=1}|(\partial\bar\partial I)(v,\bar v)|\}$, we obtain
\begin{equation}\label{eq:embedded-disc-laplacian-bound}
 |\Delta_{\rm eucl}v|\le L_I\rho_g.
\end{equation}

Let $p$ be the real-valued smooth function on $\overline{\Disc}$ such that
\begin{equation}\label{eq:disc-green-potential-problem}
 -\Delta_{\rm eucl}p=2\rho_g\quad\text{in }\Disc,\qquad p|_\T=0.
\end{equation}
Then the Green--Riesz formula \cite[pp.~356--357]{Riesz30}, \cite[Ch.~I, \S4.B.1, (4.5)]{DemCADG} gives
\begin{equation}\label{eq:disc-green-potential-representation}
 p(z)=-\int_\Disc\mathcal G_\Disc(z,w)\,d\mu_g(w).
\end{equation}
Applying the same formula componentwise to $v$ and using
\eqref{eq:embedded-disc-laplacian-bound}, we  obtain
\begin{equation}\label{eq:disc-zero-trace-green-bound}
 |v(z)|\le L_I\int_\Disc(-\mathcal G_\Disc(z,w))\rho_g(w)\,d\lambda(w)
         =\frac{L_I}{2}p(z).
\end{equation}
By \eqref{eq:green-energy},
\begin{equation}\label{eq:green-potential-energy-pairing}
 Q_h(g)=2\int_\Disc p\,d\mu_g=4\int_\Disc \rho_g p\,d\lambda.
\end{equation}
After integrating by parts, we use
\eqref{eq:embedded-disc-laplacian-bound},
\eqref{eq:disc-zero-trace-green-bound} and
\eqref{eq:green-potential-energy-pairing} to obtain
\begin{equation}\label{eq:disc-zero-trace-energy-bound}
 \|dv\|_2^2=-\int_\Disc\langle v,\Delta_{\rm eucl}v\rangle\,d\lambda
 \le\frac{L_I^2}{2}\int_\Disc \rho_g p\,d\lambda
 =\frac{L_I^2}{8}Q_h(g).
\end{equation}
Together with \eqref{eq:disc-harmonic-energy-decomposition}, we get
\begin{equation}\label{eq:closed-disc-area-green-bound}
 a_I A_h(g)\le\mathcal E(g)+\frac{L_I^2}{8}Q_h(g).
\end{equation}
Thus the first bound in \eqref{g:eq:scalar-green} holds with
$C:=\max\{1,a_I^{-1},L_I^2/(8a_I)\}$.
\end{proof}

Now we can give the proof of \Cref{g:th:main}.
\begin{proof}[Proof of \Cref{g:th:main}]
The first inequality of \eqref{g:eq:area-boundary} follows from  Dirichlet's principle
\cite[\S8.2.3, remark after Theorem~4, pp.~474--475]{Evans10}
and \eqref{eq:embedded-disc-energy-density}. 

For the second inequality of \eqref{g:eq:area-boundary}, write
$\rho_g(z):=|g'(z)|_h^2$, $\chi(z):=1-|z|^2,\ z\in\Disc$, and choose $C\ge1$ so that
\eqref{g:eq:scalar-green} gives 
\begin{equation}
A_h(g)\le C(\mathcal E(g)+\|\chi \rho_g\|_2^2).
\end{equation}
Put $\eta:=C^{-1/2}$ and $S:=\{\chi\sqrt{\rho_g}>\eta\}$.
On $\Disc\setminus S$ we have $C\chi^2 \rho_g^2\le \rho_g$. On $S$,
\eqref{g:eq:brody-speed} gives $\chi^2\rho_g^2\le D_h \rho_g$. Since
$\int_\Disc \rho_g\,d\lambda=A_h(g)/2$, it follows that
\begin{equation}\label{eq:area-on-fast-set}
 A_h(g)\le 2C\mathcal E(g)+2C D_h\int_S \rho_g\,d\lambda .
\end{equation}
By radial contraction and part~\textup{(i)} of \Cref{lem:radial-traces},
it suffices to prove the upper bound for all closed holomorphic discs
with a constant depending only on $X,h,I$. Suppose, to the contrary, that no such uniform constant exists.
Then for each integer $j\ge1,$ there is a closed holomorphic disc
$g_j\colon\Disc\to X$ such that $A_h(g_j)>j\,\mathcal E(g_j)$.
Set $\rho_{g_j}(z):=|g_j'(z)|_h^2$ and $S_j:=\{\chi\sqrt{\rho_{g_j}}>\eta\}$.
By \eqref{eq:area-on-fast-set} and the choice of $g_j$, we obtain
\begin{equation}\label{eq:fast-set-area-lower-bound}
 4CD_h\int_{S_j}\rho_{g_j}\,d\lambda\ge A_h(g_j)
 \qquad (j\ge4C).
\end{equation}
In the following, we always assume $j\geq 4C$.

For an arc $J\subseteq\T$, let $\ell_J:=\sigma(J)$ be its normalized length, and set 
$V(J):=\{s\zeta:0<s<1,\ \zeta\in J\}$. For $u\in H^{1/2}(\T;\R^N)$ set
\[
 H_J(u):=2\pi\int_J\int_\T
 \frac{|u(\zeta)-u(\xi)|^2}{|\zeta-\xi|^2}\,d\sigma(\xi)\,d\sigma(\zeta),
\]
and put
\begin{equation}\label{jdefinition}
 \tau_j(J):=\int_{V(J)}g_j^*h+jH_J(\gamma_{g_j}),\ \qquad K:=16\pi CD_h^2,\qquad m:=\lceil A_h(g_j)/K\rceil,
\end{equation}
where $\lceil x\rceil$ denotes the smallest integer greater than or equal to $x$.
Since $V(\T)=\Disc\setminus\{0\}$, Douglas's formula
\eqref{g:eq:douglas} and the choice of $g_j$ give
\begin{equation}\label{eq:total-arc-weight}
 \tau_j(\T)=A_h(g_j)+j\mathcal E(g_j)<2A_h(g_j).
\end{equation}
For $0\leq a<b\leq 2\pi$, set $J_{a,b}:=\{e^{i\theta}:a\le\theta<b\}$, then we have
\begin{equation}\label{eq:arc-weight-angular-density}
 \tau_j(J_{a,b})=\int_a^b\left(
  2\int_0^1 \rho_{g_j}(se^{i\theta})\,s\,ds
  +j\int_\T
    \frac{|\gamma_{g_j}(e^{i\theta})-\gamma_{g_j}(\xi)|^2}
         {|e^{i\theta}-\xi|^2}\,d\sigma(\xi)
 \right)d\theta .
\end{equation}
By \eqref{eq:total-arc-weight}, for fixed $a$, $\tau_j(J_{a,b})$ is a continuous increasing function of $b$.
Thus we can find 
$0=t_0<t_1<\cdots<t_m=2\pi$ such that
\begin{equation}\label{eq:arc-partition-points}
 \tau_j( J_{0,t_k})
 =\frac{k}{m}\tau_j(\T),\qquad 0\le k\le m.
\end{equation}
Here $J_{0,0}=\varnothing$.

The arcs $J_{t_{k-1},t_k}$, $1\le k\le m$, form a partition
$\mathcal J_j$ of $\T$.  By \eqref{eq:total-arc-weight} and \eqref{eq:arc-partition-points},
 for every $J\in\mathcal J_j$,
\begin{equation}\label{eq:equal-arc-weight}
 \tau_j(J)=\frac{\tau_j(\T)}{m}
 \le\frac{2A_h(g_j)}{m}\le 2K.
\end{equation}
Then \eqref{jdefinition} implies
\begin{equation}\label{eq:arc-area-energy-bounds}
 \int_{V(J)}g_j^*h\le 2K,\qquad
 H_J(\gamma_{g_j})\le\frac{2K}{j}
 \qquad (J\in\mathcal J_j).
\end{equation}

By \eqref{g:eq:brody-speed}, $\rho_{g_j}(s\zeta)\le D_h/(1-s)^2$.
Integrating in polar coordinates, we obtain
\begin{equation}\label{eq:sector-area-sum}
 \sum_{J\in\mathcal J_j}\int_{V(J)\cap\{1-|z|>\ell_J\}}\rho_{g_j}\,d\lambda
 \le 2\pi D_h\sum_{J\in\mathcal J_j}\ell_J\int_{\ell_J}^1 t^{-2}\,dt
 =2\pi D_h(m-1)<\frac{A_h(g_j)}{8CD_h}.
\end{equation}
By \eqref{eq:fast-set-area-lower-bound} and \eqref{eq:sector-area-sum},
the set
\[
 S_j\setminus\bigcup_{J\in\mathcal J_j}
 \bigl(V(J)\cap\{1-|z|>\ell_J\}\bigr)
\]
has positive $\rho_{g_j}\,d\lambda$-measure for $j\ge4C$.
Then we can choose
$J_j\in\mathcal J_j$ and $z_j=r_jp_j\in S_j\cap V(J_j)$ such that
$0<r_j<1$, $p_j$ is in the interior of $J_j$, and $1-r_j\le\ell_{J_j}$.

Set 
$$\varphi_j(w):=\frac{p_j(r_j+w)}{1+r_jw},\ w\in\Disc\qquad  f_j:=g_j\circ\varphi_j.$$
Since $\varphi_j(0)=z_j\in S_j$, the chain rule gives
\begin{equation}
|f_j'(0)|_h=(1-r_j^2)|g_j'(z_j)|_h>\eta.
\end{equation}
In at least one direction along $J_j$, the angular distance from $p_j$
to the endpoint is at least $\pi\ell_{J_j}$. Passing to a subsequence, we may assume that the chosen direction
is the same for all remaining indices $j$. Let
$E=\{e^{i\theta}:0<\theta<1/8\}$ if the direction is counterclockwise,
and take its conjugate otherwise. Since $\varphi_j(1)=p_j$ and $\varphi_j$ preserves orientation on $\T$,
this choice makes $\varphi_j(E)$ extend from $p_j$ in the direction
in which $J_j$ extends for at least $\pi\ell_{J_j}$ radians.
Put $U:=\{s\zeta:3/4<s<1,\ \zeta\in E\}$.
For $w\in U\cup E$ we have $\operatorname{Re}w>0$ and $|1-w|<1/2$.
Writing $v:=\varphi_j(w)/p_j$, we obtain
\begin{equation}\label{eq:sector-argument-bound}
 |1-v|<\frac{1-r_j}{2},\qquad
 |\arg v|\le\frac{|\operatorname{Im}v|}{\operatorname{Re}v}<1-r_j,\qquad \operatorname{Im}v=\frac{(1-r_j^2)\operatorname{Im}w}{|1+r_jw|^2},
\end{equation}
where we write $\operatorname{Re}z$ and $\operatorname{Im}z$ for the real
and imaginary parts of $z$ respectively, and use the principal
argument $\arg z\in(-\pi,\pi]$ for $z\ne0$.
Since $1-r_j\le\ell_{J_j}$, \eqref{eq:sector-argument-bound} gives
\begin{equation}
\varphi_j(U)\subset V(J_j),\qquad \varphi_j(E)\subset J_j.
\end{equation}

Put $F_j:=I\circ f_j$ and $B:=2b_IK$. Applying \eqref{eq:embedded-disc-energy-density} and
\eqref{eq:arc-area-energy-bounds}, we obtain
\begin{equation}\label{fjbound}
\int_U|dF_j|^2\,d\lambda\le b_I\int_{\varphi_j(U)}g_j^*h\le B.
\end{equation}
For a disc automorphism $\varphi$ and $\zeta,\xi\in\T$,
\begin{equation}
 |\varphi(\zeta)-\varphi(\xi)|^2
 =|\varphi'(\zeta)|\,|\varphi'(\xi)|\,|\zeta-\xi|^2 .
\end{equation}
By the fact
$\gamma_{f_j}=\gamma_{g_j}\circ\varphi_j$ and \eqref{eq:arc-area-energy-bounds},
\begin{equation}
H_E(\gamma_{f_j})=H_{\varphi_j(E)}(\gamma_{g_j})
\le H_{J_j}(\gamma_{g_j})\le2K/j.
\end{equation}
For $c_j:=\av_\T\gamma_{f_j}$, convexity of $|\cdot|^2$ and $|\zeta-\xi|\le2$ give
\begin{equation}\label{eq:arc-mean-concentration}
 \int_E|\gamma_{f_j}-c_j|^2\,d\sigma
 \le\frac{2}{\pi}H_E(\gamma_{f_j})\leq \frac{4K}{j\pi}\longrightarrow0\text{ as }j\to \infty.
\end{equation}

By \eqref{g:eq:brody-speed} and compactness of $X$, a subsequence of
$f_j$ converges locally smoothly to a holomorphic disc $f$ with $|f'(0)|_h\ge\eta$.
Since $I(X)$ is compact,  $c_j$ is a bounded sequence. By passing to a
further subsequence, we may assume $f_j$ converges locally smoothly to $f$ and  $c_j\to c$ as $j\to\infty$. Exhausting $U$ by relatively compact subsets in \eqref{fjbound}, we find that the map $F:=I\circ f$ satisfies 
\begin{equation}
\int_U|dF|^2\,d\lambda\le B.
\end{equation}
For $3/4<r<1$, \eqref{eq:radial-trace} gives 
\begin{equation}
\|F_j(r\,\cdot)-\gamma_{f_j}\|_{L^2(E)}^2\le \frac{B\log(1/r)}{2\pi}.
\end{equation}
By \eqref{eq:arc-mean-concentration}, $\gamma_{f_j}\to c$ in $L^2(E)$ as $j\to\infty$.
Letting $j\to\infty$ with $r$ fixed and then letting $r\uparrow1$, we obtain
$F(r\,\cdot)\to c$ in $L^2(E)$ as $r\uparrow1$. Applying
\eqref{eq:radial-trace} again, then we get
\begin{equation}\label{eq:radial-trace-decay}
 \int_E|F(r\zeta)-c|^2\,d\sigma(\zeta)
 \le\frac{\log(1/r)}{2\pi}\int_{U\cap\{|w|>r\}}|dF(w)|^2\,d\lambda(w)
 =o(1-r).
\end{equation}

For any closed subarc $E'\Subset E$, set
\[
 M_r:=\max_{\zeta\in E'}|F(r\zeta)-c|.
\]
 Choose $e^{i\theta_r}\in E'$ where the maximum is attained.
By \eqref{eq:embedded-disc-energy-density} and \eqref{g:eq:brody-speed}, we may take $L:=\sqrt{2b_ID_h}$ in
\begin{equation}
|\partial_\theta F(re^{i\theta})|\le \frac{L}{1-r}.
\end{equation}
Since $F$ is bounded and $E'$ has positive angular distance from
the endpoints of $E$, the arc
\[
 \left\{e^{i\theta}:
 |\theta-\theta_r|<\frac{M_r(1-r)}{2L}\right\}
\]
lies in $E$ for all $r$ sufficiently close to one. On this arc,
by the fundamental theorem of calculus, 
\begin{equation}
 |F(re^{i\theta})-c|
 \ge M_r-\frac{L}{1-r}|\theta-\theta_r|
 \ge \frac{M_r}{2}.
\end{equation}
Integrating over this arc and using \eqref{eq:radial-trace-decay},
we obtain
\begin{equation}
 \frac{M_r^3}{8\pi L}
 \le \frac{1}{1-r}\int_E|F(r\zeta)-c|^2\,d\sigma
 \longrightarrow 0\quad\text{ as }\quad r\uparrow 1.
\end{equation}
Hence $F(r\,\cdot)\to c$ uniformly on every closed subarc $E'\Subset E$.
Since $I(X)$ is closed, $c=I(x)$ for some $x\in X$.
Choose a local holomorphic chart $(V,\psi)$ around $x$ with $\psi(x)=0$.  
Since $I$ is an embedding, the above  uniform convergence implies that  $\psi\circ f$ extends
continuously to a closed subarc $E'\subset E$ with value zero. Reflect each coordinate $u$ across the arc by $\overline{u(1/\bar z)}$. The identity theorem then gives
$f\equiv x$ on $\Disc$, contradicting $|f'(0)|_h\ge\eta>0$.
\end{proof}

\subsection{Pointwise volume comparison}\label{sec:variation}

In this subsection, we prove the following pointwise volume comparison,
which we use to bound $\vol(K_X)$ from below in the proof of
\Cref{thm:main}.

\begin{theorem}\label[theorem]{v:t:main}
Let $(X,h)$ be a compact hyperbolic Hermitian manifold of dimension $n$.
If a Hermitian metric $\beta$ on $X$ satisfies
$\tr_\beta\Ric^{(1)}_\beta\ge-n$, then
\begin{equation}\label{v:out}
 \beta^n\ge \left(\frac{n+1}{D_h}\right)^n h^n ,
\end{equation}
where $D_h$ is defined in \Cref{le:brody-derivative-bound}.
\end{theorem}

We first give the following integral estimate.

\begin{lemma}\label[lemma]{le:weighted-area}
Let $g$ be a finite-area holomorphic disc in a Hermitian manifold $(X,h)$,
and suppose $\mathcal B_h(g):=\sup_{z\in\Disc}(1-|z|^2)^2|g'(z)|_h^2<\infty$.
For every biholomorphic map $\varphi\colon\Disc\to r\Disc$, $0<r\le1$,
\begin{equation}\label{v:green-area}
 \int_\Disc \log|z|^{-2}\,|(g\circ\varphi)'|_h^2\,d\lambda
 \le \sqrt{2\pi\,\mathcal B_h(g)\,A_h(g)}.
\end{equation}
\end{lemma}
\begin{proof}
Put $\rho_{g\circ\varphi}(z):=|(g\circ\varphi)'(z)|_h^2$. Since $\varphi/r$ is a disc automorphism, $(1-|z|^2)|\varphi'|=r-|\varphi|^2/r\le1-|\varphi|^2$. Hence
\begin{equation}
(1-|z|^2)^2\rho_{g\circ\varphi}(z)\le \mathcal B_h(g).
\end{equation} 
Since $2\int_\Disc \rho_{g\circ\varphi}\,d\lambda\le A_h(g)$,
we have 
\begin{equation}
\|(1-|z|^2)\rho_{g\circ\varphi}\|_2^2\le \mathcal B_h(g) A_h(g)/2.
\end{equation}
In polar coordinates, with $t=|z|^2$, expansion of $(1-t)^{-2}$ gives
$\left\|\frac{\log|z|^{-2}}{1-|z|^2}\right\|_2^2
 =\pi\int_0^1\frac{\log^2 t}{(1-t)^2}\,dt
 =2\pi\sum_{k\ge1}k^{-2}\le4\pi$.
Cauchy--Schwarz therefore gives
\begin{equation}
\int_{\mathbb D}\log|z|^{-2}\,\rho_{g\circ\varphi}\,d\lambda
\le
\left\|\frac{\log|z|^{-2}}{1-|z|^2}\right\|_2
\sqrt{\frac{\mathcal B_h(g)A_h(g)}2}
 \le \sqrt{2\pi\,\mathcal B_h(g)\,A_h(g)}.
\end{equation}
\end{proof}

\begin{lemma}\label[lemma]{v:le:radial-trace}
Let $g$ be a finite-area holomorphic disc in a compact Hermitian
manifold $(X,h)$. Then $Q_h(g)<\infty$. If
$\mathcal E(g)-\mathcal E(g(r\,\cdot))\le L(1-r)$ for some $L\ge0$
and all $r<1$ sufficiently close to $1$, then $g$ extends continuously
to $\overline\Disc$ with $\gamma_g\in W^{1,2}(\T;\R^N)$.
\end{lemma}

\begin{proof}
We may assume $g$ is nonconstant.
Put $\rho_g(z):=|g'(z)|_h^2$ and $c:=4\max\{0,\max_{|v|_h=1}H_h(v)\}$. The Chern curvature
formula for the holomorphic section $g'$ of $g^*T_X$
\cite[proof of Lemma~2.2]{BR03} gives $\Delta_{\rm eucl}\rho_g\ge-c\rho_g^2$.
Set $\varepsilon_0:=\pi/(16(1+c))$ and $C:=8/\pi$.
The scalar mean value inequality \cite[Appendix~A, Lemma~A.1 and Remarks~A.2--A.3]{RS01}
gives, for $B(a,s)\Subset\Disc$ with $\int_{B(a,s)}\rho_g\,d\lambda<\varepsilon_0$,
\begin{equation}\label{smallenergy}
s^2\rho_g(a)\le C\int_{B(a,s)}\rho_g\,d\lambda.
\end{equation}
Since $g$ has finite area, $\rho_g\in L^1(\Disc)$.
For any $a\in\Disc$, take $s:=(1-|a|)/2$. Since $B(a,s)\subset\{z\in\Disc:|z|>1-3s\}$, 
\begin{equation}
 \int_{B(a,s)}\rho_g\,d\lambda
 \le \int_{\{1-3s<|z|<1\}}\rho_g\,d\lambda \longrightarrow 0
 \qquad (|a|\to1).
\end{equation}
Since $1-|a|^2\le 4s$, \eqref{smallenergy} gives
\begin{equation}\label{estimate1}
(1-|a|^2)^2 \rho_g(a)\le 16C\int_{B(a,s)}\rho_g\,d\lambda
\le 16C\varepsilon_0
\end{equation}
whenever $|a|$ is sufficiently close to one. 
Then the continuity of $\rho_g$ on $\Disc$ implies $\mathcal B_h(g)<\infty$.
By \eqref{eq:green-energy-weighted-bound},
\begin{equation}\label{eq:finite-area-green-energy}
 Q_h(g)\le 8\int_{\mathbb D}(1-|a|^2)^2 \rho_g(a)^2\,d\lambda(a)
 \le4\mathcal B_h(g)A_h(g)<\infty.
\end{equation}

Put $g_r(z):=g(rz),z\in\Disc,$ and assume $\mathcal E(g)-\mathcal E(g_r)\le L(1-r)$
near $r=1$. For $0<r<1$, let $w\in C^\infty(\overline{r\Disc})$ solve the Dirichlet problem $-\Delta_{\rm eucl}w=c\rho_g$ in $r\Disc$ 
and $w|_{r\T}=0$. Given $y\in r\Disc$, choose a biholomorphic map
$\varphi\colon\Disc\to r\Disc$ with $\varphi(0)=y$. Apply
\eqref{eq:disc-jensen} to $w\circ\varphi$ and use \eqref{v:green-area} to obtain
\begin{equation}\label{eq:radial-poisson-potential-bound}
 0\le w(y)=\frac{c}{4\pi}\int_\Disc\log|z|^{-2}\,|(g\circ\varphi)'|_h^2\,d\lambda
 \le M:=\frac{c}{4\pi}\sqrt{2\pi\,\mathcal B_h(g)\,A_h(g)}.
\end{equation}
For $u=\rho_g\exp(-w),$ the same curvature formula \cite[proof of Lemma~2.2]{BR03} gives
\begin{equation}\label{eq:weighted-speed-subharmonic}
 \exp(w)u_{z\bar z}=|\nabla^h_{\partial_z}g'-w_zg'|_h^2
 -h(R^h(g',\overline{g'})g',g')+\tfrac{c}{4}\rho_g^2\ge0,
\end{equation}
then $u$ is subharmonic and the  circular means of $u$ are
nondecreasing.
Put $m(s):=\av_\T \rho_g(s\,\cdot)$. Since $0\le w\le M$ on $r\Disc$ and $w=0$ on $r\T$,
\begin{equation}\label{v:radial-means}
 m(s)\le \exp(M) m(r)\qquad(0<s\le r<1).
\end{equation}

Write $F:=I\circ g$, $F_r(z):=F(rz)$, $V_r(\zeta):=\partial_r F(r\zeta)$
and $E(r):=\mathcal E(g_r)=H(F_r|_\T)$. Differentiating this energy and
using harmonic orthogonality gives
$E'(r)=2\int_\Disc\langle dF_r,d\mathsf P V_r\rangle\,d\lambda$.
Since the outward normal derivative of $F_r$ on $\T$ equals $rV_r$,
integration by parts gives
\begin{equation}\label{v:eq:radial-derivative}
 E'(r)=4\pi r\|V_r\|_{L^2(\T)}^2
 -2\int_{r\Disc}\left\langle\Delta_{\rm eucl}F(z),
                    (\mathsf P V_r)(z/r)\right\rangle\,d\lambda(z).
\end{equation}
By \eqref{eq:embedded-disc-laplacian-bound} and the bounds on $dI$ in
\eqref{eq:embedded-disc-energy-density}, with $C:=\max\{1,L_I,a_I^{-1},b_I\}$ we have 
\begin{equation}\label{upper}
|\Delta_{\rm eucl}F|\le C \rho_g,\qquad C^{-1}m(r)\le\|V_r\|_{L^2(\T)}^2\le Cm(r).
\end{equation}
For $|z|=s<r$, the Poisson kernel at $z/r$ has mean one and maximum
$(r+s)/(r-s)$, which also bounds its squared $L^2(\T)$ norm. Thus \eqref{upper} gives
\begin{equation}\label{upper1}
|(\mathsf P V_r)(z/r)|\le(2Cm(r)/(r-s))^{1/2}.
\end{equation}
Applying \eqref{upper} and \eqref{upper1} to
\eqref{v:eq:radial-derivative} and integrating in polar coordinates,
we may take $c_0:=(8\pi C^2)^{-1}$ for $1/2<r<1$ to obtain
\begin{equation}\label{eq:penalty-derivative-bound}
 E'(r)\ge c_0m(r)-c_0^{-1}\sqrt{m(r)}
                     \int_0^r\frac{s\,m(s)}{\sqrt{r-s}}\,ds.
\end{equation}

Set $T_\rho:=\int_\rho^1s\,m(s)\,ds$. Since
$4\pi\int_0^1s\,m(s)\,ds=A_h(g)<\infty$, we have $T_\rho\to0$ as
$\rho\uparrow1$.
For $r>(1+\rho)/2$,
\begin{equation}\label{equa}
c_0^{-1}\int_0^\rho \frac{s\,m(s)}{\sqrt{r-s}}\,ds
\le \frac{1}{c_0}\sqrt{\frac{2}{1-\rho}}
      \int_0^\rho s\,m(s)\,ds
\le \frac{A_h(g)}{4\pi c_0}\sqrt{\frac{2}{1-\rho}}=:C_\rho.
\end{equation}
By \eqref{v:radial-means},  for any $0<\rho<r,\ b>0$,
\begin{align*}
 \int_\rho^r\frac{s\,m(s)}{\sqrt{r-s}}\,ds
 &=\int_{\substack{\rho<s<r\\ r-s<b}}
       \frac{s\,m(s)}{\sqrt{r-s}}\,ds
   +\int_{\substack{\rho<s<r\\ r-s\ge b}}
       \frac{s\,m(s)}{\sqrt{r-s}}\,ds\le \exp(M)m(r)\int_0^b\frac{dt}{\sqrt t}
      +\frac{1}{\sqrt b}\int_\rho^1 s\,m(s)\,ds\\
 &=2\exp(M)m(r)\sqrt b+\frac{T_\rho}{\sqrt b}.
\end{align*}
Hence
\begin{align*}
 \int_\rho^r\frac{s\,m(s)}{\sqrt{r-s}}\,ds&\le \inf_{b>0}\left(2\exp(M)m(r)\sqrt b
                      +\frac{T_\rho}{\sqrt b}\right)=2\sqrt{2\exp(M)m(r)T_\rho}
 \le 3\sqrt{\exp(M)m(r)T_\rho}.
\end{align*}
Combining this estimate with \eqref{eq:penalty-derivative-bound} and \eqref{equa}
gives 
\begin{equation}\label{eq:radial-energy-derivative-tail}
 E'(r)\ge \bigl(c_0-3c_0^{-1}\sqrt{\exp(M)T_\rho}\bigr)m(r)
          -C_\rho\sqrt{m(r)}.
\end{equation}
Choose $\rho$ close enough to one that
$3c_0^{-1}\sqrt{\exp(M)T_\rho}\le c_0/2$. Since
$C_\rho\sqrt{m(r)}\le (c_0/4)m(r)+C_\rho^2/c_0$, we obtain
\begin{equation}\label{derivative inequality}
E'(r)\ge c_1m(r)-C_1 \quad \text{ for }\quad (1+\rho)/2<r<1,
\end{equation} 
where $c_1:=c_0/4$ and $C_1:=C_\rho^2/c_0$.

For $r<1$ sufficiently close to $1$ and any $R\in(r,1)$, integrating
\eqref{derivative inequality} from $r$ to $R$ gives
\begin{equation}
 c_1\int_r^R m(s)\,ds
 \le \int_r^R \bigl(E'(s)+C_1\bigr)\,ds
 = E(R)-E(r)+C_1(R-r).
\end{equation}
By \Cref{lem:radial-approximation}, $E(R)\to\mathcal E(g)$ as $R\uparrow1$.
Passing to the limit and using the assumed bound
$\mathcal E(g)-E(r)\le L(1-r)$, we obtain
\begin{equation}\label{eq:integrated-radial-speed-bound}
 c_1\int_r^1 m(s)\,ds
 \le \mathcal E(g)-E(r)+C_1(1-r)
 \le (L+C_1)(1-r).
\end{equation}
Averaging \eqref{v:radial-means} over $r<s<1$ now gives
\begin{equation}\label{eq:uniform-radial-mean}
 m(r)\le\frac{\exp(M)}{1-r}\int_r^1m(s)\,ds\le \frac{\exp(M)(L+C_1)}{c_1} .
\end{equation}
Since $F(r\T)\subset I(X)$ and $|\partial_\theta F(r\zeta)|^2\le Cr^2\rho_g(r\zeta)$ by enlarging $C$ if necessary,
then \eqref{eq:uniform-radial-mean} implies  $F(r\,\cdot)$ is uniformly bounded in $W^{1,2}(\T;\R^N)$
near $r=1$. By \Cref{lem:radial-approximation}, its $L^2$-limit
is $\gamma_g$, so weak compactness gives $\gamma_g\in W^{1,2}(\T;\R^N)$.
The uniform $W^{1,2}$ bound gives equicontinuity on $\T$, so the $L^2$ convergence is uniform. 
Thus $I^{-1}\circ\gamma_g$ gives a continuous boundary extension of $g$.
\end{proof}

Now we return to the case $X$ is a compact hyperbolic manifold of dimension $n$.
Fix two Hermitian metrics $h,\beta$ on $X$, and choose a constant $\Lambda>D_h$ by \Cref{le:brody-derivative-bound}.
For a finite-area holomorphic disc $g$ in $X$, set
\[
 q_\Lambda(g):=-(n+1)\log\bigl(1-|g'(0)|_h^2/\Lambda\bigr),\qquad
 \mathcal P(g):=\frac1{2\pi}\int_\Disc\log|z|^{-2}\,g^*\Ric^{(1)}_\beta .
\]
\begin{lemma}\label[lemma]{v:l:radial}
For every $x_0\in X$ and any $0<\delta\le1/2$, the supremum
\begin{equation}\label{eq:ambient-disc-functional}
 b(\delta):=\sup_g
       \bigl(q_\Lambda(g)+\mathcal P(g)-\delta\mathcal E(g)\bigr)
\end{equation}
over all finite-area holomorphic discs $g$ centered at $x_0$ is attained.
Every maximizing holomorphic disc $g$  extends continuously to
$\overline\Disc$, and its boundary trace $\gamma_g$ belongs to $W^{1,2}(\T;\R^N)$.
\end{lemma}

\begin{proof}
By \eqref{g:eq:brody-speed}, $\mathcal B_h(g)\le D_h$.
With $M_\beta:=\max_{|v|_h=1}|\Ric^{(1)}_\beta[v]|$, we have $-M_\beta h\le\Ric^{(1)}_\beta\le M_\beta h$ on $X$.
Taking $\varphi$ to be the identity map in \eqref{v:green-area} gives
\begin{equation}\label{upper2}
|\mathcal P(g)|\le M_\beta\sqrt{2D_h A_h(g)/\pi}.
\end{equation}
By  \Cref{g:th:main}, we may find a constant 
$C_0:=C_0(X,I,h)\geq 1$ such that  $A_h(g)\le C_0\mathcal E(g)$. 
Combining this inequality with \eqref{upper2} and $|g'(0)|_h^2\le D_h$ gives
\begin{equation}\label{eq:maximizing-area-bound}
 q_\Lambda(g)+\mathcal P(g)-\delta\mathcal E(g)
 \le -(n+1)\log(1-D_h/\Lambda)
     +M_\beta\sqrt{\frac{2D_h A_h(g)}{\pi}}
     -\frac{\delta}{C_0}A_h(g).
\end{equation}
For fixed $0<\delta\leq 1/2$ and $\Lambda>D_h$, the right-hand side is bounded
above and tends to $-\infty$ as $A_h(g)\to\infty$.
By considering  the constant disc $g\equiv x_0$, $0\le b(\delta)<\infty$.
Choose a maximizing sequence $(g_j)$ of holomorphic discs in $X$ centered at $x_0$ such that
\begin{equation}
 q_\Lambda(g_j)+\mathcal P(g_j)-\delta\mathcal E(g_j)
 \ge b(\delta)-1/j\ge-1\qquad (j\geq 1).
\end{equation}
Then \eqref{eq:maximizing-area-bound} gives $\sup_j A_h(g_j)<\infty$.
Since $X$ is compact, \eqref{g:eq:brody-speed} gives a subsequence converging locally uniformly,
hence locally smoothly in holomorphic charts, to a holomorphic disc $g$ centered at $x_0$. Thus
$q_\Lambda(g_j)\to q_\Lambda(g)$. By \Cref{lem:bounded-area-limits},
$g$ has finite area,
\begin{equation}
\mathcal P(g_j)\to\mathcal P(g),\qquad \mathcal E(g)\le\liminf_{j\rightarrow\infty}\mathcal E(g_j).
\end{equation} 
Hence  $b(\delta)=q_\Lambda(g)+\mathcal P(g)-\delta\mathcal E(g)$.

Now let $g$ be any maximizing disc, and set $g_r(z):=g(rz)$ for $1/2\le r<1$.
Changing of variables gives
\begin{equation}\label{eq:radial-curvature-integral}
 \mathcal P(g_r)=\frac1{2\pi}\int_{r\Disc}
       \log\frac{r^2}{|w|^2}\,g^*\Ric^{(1)}_\beta .
\end{equation}
Extend the weight $W_r(w):=\log(r^2/|w|^2)$ in
\eqref{eq:radial-curvature-integral} by zero on $\Disc\setminus r\Disc$.
Since $r\ge1/2$, for $0<|w|<1,$ we have
\begin{equation}
 0\le \log|w|^{-2}-W_r(w)
 =-2\log\max\{r,|w|\}
 \le -2\log r
 =2\int_r^1\frac{dt}{t}\le 4(1-r).
\end{equation}
Put $d:=|g'(0)|_h^2<\Lambda$. Using
$-M_\beta h\le\Ric^{(1)}_\beta\le M_\beta h$ and $g_r'(0)=rg'(0)$, we obtain
\begin{equation}\label{eq:radial-weights}
 \begin{aligned}
 |\mathcal P(g)-\mathcal P(g_r)|
 &\le \frac{M_\beta A_h(g)}{\pi}\log\frac1r
 \le \frac{2M_\beta A_h(g)}{\pi}(1-r),\\
 0\le q_\Lambda(g)-q_\Lambda(g_r)
 &=(n+1)\log\left(1+\frac{(1-r^2)d}{\Lambda-d}\right)
 \le \frac{2(n+1)d}{\Lambda-d}(1-r).
 \end{aligned}
\end{equation}
By maximality of $g$,
\begin{equation}\label{eq:radial-maximizer-comparison}
 \delta\bigl(\mathcal E(g)-\mathcal E(g_r)\bigr)
 \le q_\Lambda(g)-q_\Lambda(g_r)+\mathcal P(g)-\mathcal P(g_r)
 \le \left(\frac{2(n+1)d}{\Lambda-d}
           +\frac{2M_\beta A_h(g)}{\pi}\right)(1-r).
\end{equation}
The boundary assertion follows from \Cref{v:le:radial-trace} as $\delta>0$ is fixed.
\end{proof}

We say a frame $B=(B_1,\ldots,B_n)$ along $g\colon\overline\Disc\to(X,\beta)$
is \emph{$\beta$-unitary on $\T$} if $B_1(\zeta),\ldots,B_n(\zeta)$ form
an orthonormal basis of $(T^{1,0}_{g(\zeta)}X,\beta_{g(\zeta)})$ for every $\zeta\in\T$.

\begin{lemma}\label[lemma]{v:l:factor}
Let $(X,\beta)$ be a compact Hermitian $n$-fold and let
$g\colon\overline\Disc\to X$ be continuous, holomorphic on $\Disc$, with
$g|_\T\in W^{1,2}(\T,X)$. Then $g\in W^{1,4}(\Disc,X)$ and $g^*T_X$
has a holomorphic frame $B=(B_1,\ldots,B_n)\in W^{1,4}(\Disc,g^*T_X)$
which is continuous on $\overline\Disc$ and $\beta$-unitary  on $\T$.
Let $x_0$ be the center of $g$. On some ball $P\subset\C^n$ about $0$
there is a family $t\mapsto g_t$ of class $C^2(P,W^{1,4}(\Disc,X))$,
holomorphic in $(z,t)\in\Disc\times P$, with $g_0=g$. Each $g_t$ is centered at $x_0$, and, for $1\le j\le n$,
\begin{equation}\label{v:eq:centre-variations}
 V_j(z):=\left.\partial_{t_j}g_t(z)\right|_{t=0}=zB_j(z),\qquad
 \left.\partial_{t_j}g_t'(0)\right|_{t=0}=B_j(0).
\end{equation}
\end{lemma}
\begin{proof}
The graph of $g$ over the closed disc is
\[
 \Gamma_g:=\{(z,g(z)):z\in\overline\Disc\}\subset\C\times X .
\]
We use graph coordinates to construct a frame along $g$, as in
\cite[Lemma~5.2 and Step~1 of the proof of Theorem~5.5]{LSWZ26}.
By Forstneri\v{c}'s graph theorem \cite[Theorem~1.2]{Forstneric},
there exist a holomorphic vector bundle $p\colon \mathcal{V}\to U_0$ over an open
neighbourhood $U_0$ of $\overline\Disc$ and a biholomorphic map  $\Phi$ from
an open neighbourhood $V\subset U_0\times X$ of $\Gamma_g$  onto an open subset of $\mathcal{V}$
such that $p(\Phi(z,y))=z$ for every $(z,y)\in V.$ Choose $\varepsilon>0$ such that $(1+\varepsilon)\Disc\subset U_0$.
By \cite[Satz~3]{Rohrl57} (see also \cite[Theorem~30.4]{Forster81}), $\mathcal{V}$ is holomorphically trivial
over this disc, so in this trivialization we may write
\[
 \Phi(z,y)=(z,\phi_z(y)),\qquad
 \Phi^{-1}(z,w)=(z,\phi_z^{-1}(w)).
\]
In these coordinates the graph of $g$ is
\[
 \Phi(\Gamma_g)=\{(z,s(z)):z\in\overline\Disc\},\qquad
 s(z):=\phi_z(g(z)).
\]
Thus, for any $z\in\overline{\Disc}$, we have a linear isomorphism
\[
 E(z):=\left.d_w(\phi_z^{-1})\right|_{w=s(z)}
 \colon\C^n\xrightarrow{\cong}T^{1,0}_{g(z)}X .
\]
Let $e_1,\ldots,e_n$ be the standard basis of $\C^n$. For
$z\in\overline\Disc$, set
\[
 E_j(z):=E(z)e_j\in T^{1,0}_{g(z)}X,\qquad 1\le j\le n .
\]
Since $E(z)$ is an isomorphism, $E_1(z),\ldots,E_n(z)$ form a basis
of $T^{1,0}_{g(z)}X$ for every $z\in\overline{\Disc}$. Thus $E_1,\ldots,E_n$ form a
frame of $g^*T_X$, which is holomorphic in $\Disc$ and continuous on
$\overline\Disc$. 

The chain rule on $\T$ gives $s|_\T\in W^{1,2}(\T;\C^n)$.
For each component $s_\ell(z)=\sum_{k\ge0}a_kz^k$ of $s$, the Parseval identity gives
\begin{equation}
 \|s_\ell'\|_{H^2}^2=\sum_{k\ge1}k^2|a_k|^2
 =\|\partial_\theta(s_\ell|_\T)\|_2^2 .
\end{equation}
Carleman's inequality \cite{Car21} gives
\begin{equation}\label{eq:hardy-area-four}
 \frac1\pi\int_\Disc |s_\ell'|^4\,d\lambda
 \le\|s_\ell'\|_{H^2}^4 .
\end{equation}
Thus $s\in W^{1,4}(\Disc;\C^n)$.  Since
$s\in W^{1,4}(\Disc;\C^n)$ and $s|_\T\in W^{1,2}(\T;\C^n)$, the
chain rule gives 
\begin{equation}
 g\in W^{1,4}(\Disc,X),\quad E_j\in W^{1,4}(\Disc,g^*T^{1,0}X),\quad
 E_j|_\T\in W^{1,2}\bigl(\T,(g|_\T)^*T_X\bigr),\quad (1\leq j\leq n).
\end{equation}

Since $g\in W^{1,4}(\Disc,X)$, \cite[Theorem~5.5]{LSWZ26} with $p=4$
and the Hermitian metric $\beta$ gives a holomorphic frame
$B=(B_1,\ldots,B_n)$ of $g^*T_X$ of class $W^{1,4}$ that is
$\beta$-unitary on $\T$. This theorem applies to any smooth Hermitian
metric, since its proof uses only the uniform positivity and
$W^{1,4}$ regularity of the Gram matrix in the initial frame.
By the Sobolev embedding theorem, $B$ is continuous on
$\overline\Disc$. Regard $B(z)\colon\C^n\to T^{1,0}_{g(z)}X$ as the map $e_j\mapsto B_j(z)$.
In the graph coordinates, set
\[
 U(z):=E(z)^{-1}B(z),\qquad z\in\overline\Disc,
\]
then $U$ is holomorphic in $\Disc$. The chain and product rules give
$U\in W^{1,4}(\Disc;\C^{n\times n})$. As $E$ and $B$ are continuous frames on $\overline\Disc$, $U$ is continuous and invertible there.

Clearly, we may choose
an open ball $P\subset\C^n$ centered at $0$ such that $(z,s(z)+zU(z)t)$ lies
in the domain of $\Phi^{-1}$ for all
$(z,t)\in\overline\Disc\times\overline P$. Define
\[
 g_t(z):=\phi_z^{-1}\bigl(s(z)+zU(z)t\bigr).
\]
Since $\phi_z^{-1}(s(z))=g(z)$, we have
\begin{equation}\label{eq:centred-family-graph}
 \phi_z(g_t(z))-s(z)=zU(z)t,\qquad g_0=g.
\end{equation}
Setting $z=0$ shows that each $g_t$ is centered at $x_0$.
The map $(z,t)\mapsto g_t(z)$ is holomorphic on $\Disc\times P$
and is holomorphic in $t$ for each $z\in\T$.
We verify that $t\mapsto g_t$ is of class $C^2$ in $W^{1,4}(\Disc,X)$.
After composition with the fixed smooth embedding $I\colon X\to\R^N$,
differentiate the formula for $I\circ g_t$ up to twice in the real
parameter variables and once in the disc variables. By the chain rule,
each resulting term contains at most one factor from $ds$ or $dU$.
These factors belong to $L^4$, and their coefficients are bounded
and depend continuously on $t$ in the uniform norm. Thus the parameter
derivatives through order two depend continuously on $t$ in $W^{1,4}$.

Finally, differentiating in $t_j$ at $t=0$ gives
\begin{equation}
 \left.\partial_{t_j}g_t(z)\right|_{t=0}
 =d_w(\phi_z^{-1})_{s(z)}\bigl(zU(z)e_j\bigr)
 =zE(z)U(z)e_j=zB_j(z).
\end{equation}
Since all the discs $g_t$ have center $x_0$, the vectors $g_t'(0)$
all lie in the same tangent space $T_{x_0}^{1,0}X$. Choose a holomorphic
chart $x$ at $x_0$ and set $F(z,t):=x(g_t(z))$ near $(0,0)$. Then
\begin{equation}
 \partial_z F(0,t)=dx_{x_0}\bigl(g_t'(0)\bigr),\qquad
 \left.\partial_{t_j}F(z,t)\right|_{t=0}
 =z\,dx_{g(z)}\bigl(B_j(z)\bigr).
\end{equation}
Since $F$ is jointly holomorphic, we may interchange the $z$ and
$t_j$ derivatives. Thus
\begin{equation}
 dx_{x_0}\!\left(\left.\partial_{t_j}g_t'(0)\right|_{t=0}\right)
 =\left.\partial_z\!\left(z\,dx_{g(z)}(B_j(z))\right)\right|_{z=0}
 =dx_{x_0}(B_j(0)).
\end{equation}
The last equality follows from the product rule at $z=0$.
Since $dx_{x_0}$ is an isomorphism, we obtain
\begin{equation}
\left.\partial_{t_j}g_t'(0)\right|_{t=0}=B_j(0).
\end{equation}
This proves \eqref{v:eq:centre-variations}.
\end{proof}

Along the family
\eqref{eq:centred-family-graph}, write
\[
 \mathcal L\Phi:=
 \left.\sum_{j=1}^n\partial_{t_j}\partial_{\bar t_j}\Phi(g_t)\right|_{t=0}.
\]

\begin{lemma}\label[lemma]{v:l:hess}
Let $h,\beta$ be two Hermitian metrics on a compact complex manifold $X$ of
dimension $n$, and let $I\colon X\hookrightarrow\R^N$ be the smooth
embedding defining $\mathcal E$. For the family $g_t$ of
\Cref{v:l:factor}, the map $t\mapsto\mathcal E(g_t)$ is $C^2$ and
there is a constant $C_H:=C_H(X,h,\beta,I)>0$ such that 
\begin{equation}\label{v:EH}
 \mathcal L\mathcal E\le C_H(1+A_h(g)).
\end{equation}
\end{lemma}
\begin{proof}
Put
\[
 M:=B^*\beta(g)B,\quad m_j:=M_{jj},\quad R_j:=\frac{2B_j}{1+m_j},\quad
 \psi:=\sum_{j=1}^n\log\frac{1+m_j}{2}-\frac12\log\det M .
\]
Here $M_{ij}(z)=\beta_{g(z)}(B_j(z),B_i(z))$. Then $|R_j|_\beta\le1$ on $\Disc$ and $\psi\geq 0$ since $(1+m_j)/2\ge\sqrt{m_j}$ and $\det M\le\prod_jm_j$.
Since $B=(B_1,\cdots,B_n)$ is $\beta$-unitary on $\T$, $R_j=B_j$ on $\T$ and $\psi|_\T=0$.

Write $\rho:=|g'|_h^2$ and $\nabla:=g^*\nabla^\beta$. For a section $S$ along $g$,
use $|\nabla S|^2:=|\nabla_{\partial_a}S|_\beta^2+|\nabla_{\partial_b}S|_\beta^2$.
With $K:=\max_{|v|_h=|w|_\beta=1}|\beta(R^\beta(v,\bar v)w,w)|$, the maximum
being over pairs in the same fibre, we have
\begin{equation}\label{curvaturebound}
|\beta(R^\beta(v,\bar v)w,w)|\le K|v|_h^2|w|_\beta^2.
\end{equation}
For a scalar function $u$, we write $u_z:=\partial_z u$ and
$u_{z\bar z}:=\partial_{\bar z}\partial_z u$.
For each $j$ put $\xi_j:=\nabla_{\partial_z}B_j$. The Chern curvature
formula \cite[proof of Lemma~2.2]{BR03} and
$|\beta(\xi_j,B_j)|^2\le m_j|\xi_j|_\beta^2$ give
\begin{equation}\label{estimatederivative}
 \begin{aligned}
 (m_j)_z&=\beta(\xi_j,B_j),\qquad
 (m_j)_{z\bar z}=|\xi_j|_\beta^2-\beta(R^\beta(g',\overline{g'})B_j,B_j),\\
 (\log(1+m_j))_{z\bar z}
 &=\frac{(1+m_j)(m_j)_{z\bar z}-|(m_j)_z|^2}{(1+m_j)^2}
 \ge\frac{|\xi_j|_\beta^2}{(1+m_j)^2}-K\rho .
 \end{aligned}
\end{equation}

For each $z\in\Disc$, define
\[
 F_z\colon T^{1,0}_{g(z)}X\longrightarrow T^{1,0}_{g(z)}X,\qquad
 F_z(w):=\frac{2w}{1+|w|_\beta^2}.
\]
For $w,v\in T^{1,0}_{g(z)}X$, its real differential satisfies
\begin{equation}
 |(dF_z)_w(v)|_\beta\le\frac{2}{1+|w|_\beta^2}\,|v|_\beta .
\end{equation}
Since  $R_j(z)=F_z(B_j(z))$ and $\nabla\beta=0$, a direct computation gives 
\begin{equation}
 |\nabla R_j|\le \frac{2}{1+m_j}|\nabla B_j|.
\end{equation}
By the holomorphicity of $B_j$, we have $|\nabla B_j|^2=2|\xi_j|_\beta^2$. Using \eqref{estimatederivative},
we obtain
\begin{equation}
 |\nabla R_j|^2
 \le \frac{8|\xi_j|_\beta^2}{(1+m_j)^2}
 \le 2\Delta_{\rm eucl}\log(1+m_j)+8K\rho .
\end{equation}

Since $B$ is a holomorphic frame, 
\begin{equation}
 \Delta_{\rm eucl}\log\det M=-4\Ric^{(1)}_\beta[g'].
\end{equation}
By \eqref{curvaturebound}, $|\Ric^{(1)}_\beta[g']|\le nK\rho$. By the
definition of $\psi$, we then find
\begin{equation}\label{negative2}
 \sum_{j=1}^n|\nabla R_j|^2
 \le 2\Delta_{\rm eucl}\psi+12nK\rho .
\end{equation}
Since $g',\xi_j\in L^4(\Disc)$ and the curvature of $\beta$ is bounded, the preceding curvature identities give $\Delta_{\rm eucl}\psi\in L^2(\Disc)$.
Since $\psi\ge0$ and $\psi|_\T=0$, \eqref{eq:boundary-flux-sign}
gives 
\begin{equation}\label{negative1}
\int_\Disc\Delta_{\rm eucl}\psi\,d\lambda\le0.
\end{equation}
By \eqref{negative1} and \eqref{negative2},
\begin{equation}\label{v:column-energy}
 \sum_{j=1}^n\int_\Disc|\nabla R_j|^2\,d\lambda
 \le 12nK\int_\Disc\rho\,d\lambda
 =6nK A_h(g).
\end{equation}
By \Cref{v:l:factor}, $t\mapsto I\circ g_t$ is $C^2$ in
$W^{1,4}(\Disc;\R^N)$. Since the trace operator
$\operatorname{Tr}_\T\colon W^{1,4}(\Disc;\R^N)\to H^{1/2}(\T;\R^N)$
is bounded and linear, $\gamma_{g_t}=\operatorname{Tr}_\T(I\circ g_t)$
depends $C^2$ on $t$.
For each $\zeta\in\T$, the map $t\mapsto g_t(\zeta)$ is holomorphic
and the vectors $V_j(\zeta)=\zeta B_j(\zeta)$ form a $\beta$-unitary
basis. Applying $\Delta_\beta$ to $I$ componentwise, the chain rule gives
\begin{equation}
 \left.\partial_{t_j}\gamma_{g_t}\right|_{t=0}
 =\partial I(V_j)|_\T,\qquad
 \left.\sum_{j=1}^n\partial_{t_j}\partial_{\bar t_j}
 \gamma_{g_t}\right|_{t=0}
 =(\Delta_\beta I)\circ(g|_\T).
\end{equation}
Differentiating $\mathcal E(g_t)=H(\gamma_{g_t})$ and applying the
Cauchy--Schwarz inequality, we obtain 
\begin{equation}\label{v:eq:trace-hessian}
 \mathcal L\mathcal E
 \le 2\sum_{j=1}^n H(\partial I(V_j)|_\T)
 +2H(\gamma_g)^{1/2}
   H((\Delta_\beta I)\circ(g|_\T))^{1/2}.
\end{equation}
By Dirichlet's principle and the Cauchy--Schwarz inequality, the last term of
\eqref{v:eq:trace-hessian} satisfies
\begin{equation}\label{eq:trace-hessian-second-term}
 2H(\gamma_g)^{1/2}H((\Delta_\beta I)\circ(g|_\T))^{1/2}\le C A_h(g)
\end{equation}
for some constant $C:=C(X,h,\beta,I)>0$.
To estimate the first term in \eqref{v:eq:trace-hessian}, set
\[
 Y_j(z):=z\,\partial I_{g(z)}(R_j(z)).
\]
Since $R_j=B_j$ on $\T$, we have
$Y_j|_\T=\partial I(V_j)|_\T$. Since  $|R_j|_\beta\le1$, $g$ is holomorphic and $X$ is compact, 
by enlarging $C$ if necessary,
\begin{equation}
 |dY_j|\le C\bigl(1+|g'|_h+|\nabla R_j|\bigr).
\end{equation}
Squaring and integrating, we obtain
\begin{equation}
 \int_\Disc |dY_j|^2\,d\lambda
 \le 12C^2\left(1+A_h(g)+\int_\Disc|\nabla R_j|^2\,d\lambda\right).
\end{equation}
Thus $Y_j\in W^{1,2}(\Disc;\C^N)$ by \eqref{v:column-energy}.
Applying Dirichlet's principle again,
\begin{equation}\label{eq:trial-extension-boundary-energy}
 \sum_{j=1}^n H(\partial I(V_j)|_\T)
 \le \sum_{j=1}^n\int_\Disc |dY_j|^2\,d\lambda
 \le 12nC^2(1+6K)(1+A_h(g)).
\end{equation}
Substituting \eqref{eq:trial-extension-boundary-energy} and
\eqref{eq:trace-hessian-second-term} into \eqref{v:eq:trace-hessian}
proves \eqref{v:EH} with $C_H:=24nC^2(1+6K)+C$.
\end{proof}

\begin{proposition}\label[proposition]{prop:relative-volume}
Let $X$ be a compact complex $n$-fold and $\eta$ a smooth semipositive
$(1,1)$-form. Assume that $\eta>0$ on $X$ or $X$ is K\"ahler, and fix
$D>0$ such that
\begin{equation}\label{eq:source-centre-bound}
 |f'(0)|_\eta^2\le D\qquad\text{for every holomorphic }f\colon\Disc\to X.
\end{equation}
Let $\beta$ be a Hermitian metric with $\tr_\beta\Ric^{(1)}_\beta\ge-n$, and
let $\mathcal P$ be the Green integral of $\Ric^{(1)}_\beta$ defined before \Cref{v:l:radial}.
For each $x_0\in X$ with $\eta(x_0)>0$, put $D_{x_0}:=\sup_q|q'(0)|_\eta^2$,
where $q$ ranges over holomorphic discs in $X$ with center $x_0$.
Then for every finite-area disc $f$ in $X$ with center $x_0$,
\begin{equation}\label{eq:fixed-disc-volume}
 \beta^n(x_0)\ge\left(\frac{n+1}{D_{x_0}}\right)^n
 e^{\mathcal P(f)}\eta^n(x_0).
\end{equation}
In particular,
\begin{equation}\label{eq:vol}
 \beta^n\ge\left(\frac{n+1}{D}\right)^n\eta^n .
\end{equation}
\end{proposition}

In the following proof $C_0$ denotes the area constant in
\eqref{g:eq:area-boundary} if $\eta>0$ everywhere, and in \eqref{g:eq:scalar-green} otherwise.

\begin{proof}
\textbf{Step 1. Existence and boundary regularity of a maximizing disc.}
Fix a smooth embedding $I$ defining $\mathcal E$. If $\eta>0$ everywhere,
take $h:=\eta$ and $\mathcal A:=\mathcal E$. Otherwise fix a K\"ahler
metric $h$ and take $\mathcal A:=\mathcal E+Q_h$. 
Then $\mathcal{A}(f)$ is finite for every finite-area holomorphic disc $f$ by \Cref{v:le:radial-trace} and the trace theorem.
It suffices to work at a point $x_0$ where $\eta(x_0)>0$, since $\eta^n$
vanishes elsewhere.
Small coordinate discs and \eqref{eq:source-centre-bound} give $0<D_{x_0}\le D$.
For $\Lambda>D_{x_0}$ and $0<\delta\le1/2$, define
\[
 \begin{aligned}
 q_\Lambda(f)&:=-(n+1)\log(1-|f'(0)|_\eta^2/\Lambda),\\
 F_\delta(f)&:=q_\Lambda(f)+\mathcal P(f)-\delta\mathcal A(f),\qquad
 b(\delta):=\sup_f F_\delta(f),
 \end{aligned}
\]
where $f$ ranges over finite-area holomorphic discs with center $x_0$.

If $h=\eta>0$, let $F\colon\C\to X$ be an entire curve. For each
$a\in\C$, applying \eqref{eq:source-centre-bound} to $z\mapsto F(a+Rz)$
gives $|F'(a)|_h^2\le D/R^2$ for all $R>0$. Letting $R\to\infty$
gives $F'(a)=0$, so $F$ is constant. Thus $X$ is hyperbolic and
\Cref{g:th:main} gives $A_h(f)\le C_0\mathcal E(f)$.
The proof of \Cref{v:l:radial} applies with $\Lambda>D_{x_0}$.
Since $f$ has center $x_0$, the bound $|f'(0)|_h^2\le D_{x_0}<\Lambda$
bounds $q_\Lambda(f)$ uniformly and suffices for the radial comparison
\eqref{eq:radial-weights}. The global bound $\mathcal B_h(f)\le D_h\le D$
still gives \eqref{upper2} and local compactness. Thus a maximizer $g$
exists, extends continuously to $\overline\Disc$, and satisfies
$g|_\T\in W^{1,2}(\T,X)$. It remains to find a maximizer with the same boundary regularity
when $\eta$ may degenerate. In this case $h$ is K\"ahler and
$\mathcal A=\mathcal E+Q_h$. By \eqref{g:eq:scalar-green},
$A_h(f)\le C_0\mathcal A(f)$, and $Q_h(f)\le\mathcal A(f)$ by definition.
To bound $\mathcal P(f)$, we first estimate $|f'|_h$ near the center.

Let $a_{x_0}$ be half the smallest eigenvalue of $\eta(x_0)$ relative to $h(x_0)$.
Choose a coordinate neighbourhood $U$ of $x_0$ with $\eta\ge a_{x_0}h$ on $U$.
For $z\in\Disc$, let $\phi_z(w):=(z+w)/(1+\bar z w)$.
Applying \eqref{eq:source-centre-bound} to $f\circ\phi_z$, we obtain
\begin{equation}\label{derivative1}
 (1-|z|^2)^2|f'(z)|_\eta^2\le D .
\end{equation}
Choose $0<r<1$ such that
\begin{equation}
 \sqrt{\frac{2D}{a_{x_0}}}\int_0^r\frac{ds}{1-s^2}
 <\dis_h(x_0,X\setminus U).
\end{equation}
For $\theta\in\mathbb R$, put $\gamma_\theta(s):=f(se^{i\theta})$.
If $\gamma_\theta$ first leaves $U$ at $t\le r$, then
\begin{equation}
 \dis_h(x_0,X\setminus U)
 \le \operatorname{Length}_h(\gamma_\theta|_{[0,t]})
 =\sqrt{2}\int_0^t |f'(se^{i\theta})|_h\,ds
 \le \sqrt{\frac{2D}{a_{x_0}}}\int_0^t\frac{ds}{1-s^2},
\end{equation}
a contradiction. Thus $f(r\Disc)\subset U$. By $\eta\ge a_{x_0}h$
on $U$ and \eqref{derivative1},
\begin{equation}
 \sup_{z\in r\Disc}|f'(z)|_h
 \le \frac{\sqrt{D/a_{x_0}}}{1-r^2}.
\end{equation}
For $0<R<1$, let $p_R\in C^\infty(\overline{R\Disc})$ solve
\begin{equation}
 -\Delta_{\rm eucl}p_R=2|f'|_h^2 \quad\text{in }R\Disc,\qquad
 p_R=0 \quad\text{on }R\T.
\end{equation}
Set $w_R(z):=R^2-|z|^2$, then integrating by parts and \eqref{enery2} give
\begin{equation}
 \int_{R\Disc} w_R\,f^*h
 =\int_{R\Disc}\langle dw_R,dp_R\rangle\,d\lambda\leq R^2\sqrt{\pi Q_h(f(R\,\cdot))}.
\end{equation}
By the Monotone convergence theorem and
\eqref{eq:contracted-green-kernel}, we obtain
\begin{equation}
\int_\Disc(1-|z|^2)f^*h\le\sqrt{\pi Q_h(f)}.
\end{equation}
Put $M_\beta:=\max_{|v|_h=1}|\Ric^{(1)}_\beta[v]|$, so that
$-M_\beta h\le\Ric^{(1)}_\beta\le M_\beta h$. Then
\begin{equation}\label{eq:intro-truncated-green-weight}
 \begin{aligned}
 |\mathcal P(f)|
 &\le \frac{M_\beta}{2\pi}\int_\Disc \log|z|^{-2}\,f^*h \\
 &\le \frac{M_\beta}{2\pi}\left(
       \frac{2D}{a_{x_0}(1-r^2)^2}
       \int_{r\Disc}\log|z|^{-2}\,d\lambda
       +\frac{1}{r^2}\int_\Disc(1-|z|^2)\,f^*h\right) \\
 &\le \frac{M_\beta}{2\pi}\left(
       \frac{2\pi D r^2(1+\log r^{-2})}{a_{x_0}(1-r^2)^2}
       +\frac{\sqrt{\pi Q_h(f)}}{r^2}\right) \\
 &\le C_{x_0}\bigl(1+\sqrt{Q_h(f)}\bigr),
 \end{aligned}
\end{equation}
where, using $r^2(1+\log r^{-2})\le1$, we take
\[
 C_{x_0}:=M_\beta\left(\frac{D}{a_{x_0}(1-r^2)^2}+\frac1{r^2}\right).
\]
Since $|f'(0)|_\eta^2\le D_{x_0}<\Lambda$ and $Q_h(f)\le\mathcal A(f)$,
\begin{equation}
 F_\delta(f)\le C_{x_0,\Lambda}
      +C_{x_0}\sqrt{\mathcal A(f)}-\delta\mathcal A(f),
\end{equation}
where 
\[
C_{x_0,\Lambda}:=(n+1)\log\frac{\Lambda}{\Lambda-D_{x_0}}+C_{x_0}.
\]
For fixed $\delta>0$ this upper bound tends to $-\infty$ as
$\mathcal A(f)\to\infty$. The constant disc $f\equiv x_0$ satisfies $F_\delta(f)=0$, so $0\le b(\delta)<\infty$.
The preceding bound implies that every maximizing sequence $(g_j)$
satisfies $\sup_j\mathcal A(g_j)<\infty$. Hence, by
\eqref{g:eq:scalar-green},
\begin{equation}
 \sup_j\bigl(A_h(g_j)+Q_h(g_j)\bigr)<\infty.
\end{equation}

Since $h$ is K\"ahler and $A_h(g_j)$ and $Q_h(g_j)$ are uniformly
bounded,  by taking a subsequence if necessary, \cite[Proposition~3.9]{LSWZ26} allows us to assume $g_j$ converging 
locally smoothly to a finite-area holomorphic disc $g$ with center $x_0$ and
$Q_h(g)\le\liminf_j Q_h(g_j)$. Since $g_j'(0)\to g'(0)$ and
$D_{x_0}<\Lambda$, we have $q_\Lambda(g_j)\to q_\Lambda(g)$.
By \Cref{lem:bounded-area-limits}, $\mathcal P(g_j)\to\mathcal P(g)$
and $\mathcal E(g)\le\liminf_j\mathcal E(g_j)$. Hence,
\begin{equation}
 b(\delta)\ge F_\delta(g)\ge\limsup_j F_\delta(g_j)=b(\delta),
\end{equation}
and then  $g$ is a maximizer.
For $1/2<R<1$, set $g_R(z):=g(Rz)$, then we have $Q_h(g_R)\le Q_h(g)$ by
\eqref{eq:contracted-green-kernel}. Maximality and \eqref{eq:radial-weights},
with $|g'(0)|_\eta^2<\Lambda$, give
\begin{equation}
 \delta\bigl(\mathcal E(g)-\mathcal E(g_R)\bigr)
 \le q_\Lambda(g)-q_\Lambda(g_R)+\mathcal P(g)-\mathcal P(g_R)
 \le \left(\frac{2(n+1)|g'(0)|_\eta^2}{\Lambda-|g'(0)|_\eta^2}
       +\frac{2M_\beta A_h(g)}{\pi}\right)(1-R).
\end{equation}
By \Cref{v:le:radial-trace}, $g$ extends continuously to $\overline\Disc$
and $g|_\T\in W^{1,2}(\T,X)$.

\textbf{Step 2. Second variation of $\mathcal A$.}
In either case let $g=g_\delta$ be a maximizing disc and take the family
$g_t$ of \Cref{v:l:factor}. For $\mathcal A=\mathcal E+Q_h$,
each family $g_{\tau e_j}$ is jointly holomorphic in $(\tau,z)$ and
$C^2$ in $\tau$ with values in $W^{1,4}(\Disc,X)$. 
For the same K\"ahler form, the Green kernel and squared tangent norm of
\cite{LSWZ26} are $4\pi\mathcal G_\Disc$ and $2|\cdot|_h^2$, respectively.
Thus $Q_h=-\Psi/(2\pi)$ for the functional $\Psi$ of \cite[Definition~3.1]{LSWZ26},
and \cite[Proposition~4.4(2), (4.4)]{LSWZ26} gives
\begin{equation}\label{eq:relative-green-hessian}
 \mathcal L Q_h
 \le 2A_h(g)\sum_{j=1}^n\max_\T |V_j|_h^2
 \le 2n\max_X(\tr_\beta h)A_h(g),
\end{equation}
where $V_j(\xi):=\xi B_j(\xi)$, and the last inequality uses $|V_j|_\beta=1$ on $\T$ (see \Cref{v:l:factor}).
Using \eqref{g:eq:area-boundary} and \eqref{v:EH} when $\mathcal A=\mathcal E$,
and \eqref{g:eq:scalar-green}, \eqref{v:EH} and \eqref{eq:relative-green-hessian}
when $\mathcal A=\mathcal E+Q_h$, we obtain
\begin{equation}\label{eq:source-area-functional-hessian}
 A_h(g)\le C_0\mathcal A(g),\qquad
 \mathcal L\mathcal A\le C_1(1+A_h(g)).
\end{equation}
In both cases we may take
$C_1:=24nC^2(1+6K)+C+2n\max_X\tr_\beta h$, with $C,K$ as in the proof of
\Cref{v:l:hess}. Both bounds are uniform in $x_0,\Lambda,\delta$ and $g$.

\textbf{Step 3. Determinant comparison and passage to the limit.}
Put $B_{x_0}:=B(0)$, $W_{x_0}:=\log((\beta^n/\eta^n)(x_0))$ and
$u:=\log\det(B^*\beta(g)B)$. Holomorphicity of $B$ gives
$2\pi\ddc u=-g^*\Ric^{(1)}_\beta$, and boundary unitarity gives $u|_\T=0$.
Since $g$ has finite $h$-area, $\Delta_{\rm eucl}u\in L^1(\Disc)$.
Hence \eqref{eq:disc-jensen} gives $u(0)=\mathcal P(g)$. 
By definition, at $x_0$, 
\begin{equation}\label{determiant}
\frac{\det(B_{x_0}^*\eta B_{x_0})}{\det(B_{x_0}^*\beta B_{x_0})}=e^{-W_{x_0}}.
\end{equation}
For $0<r<1$, apply \cite[Proposition~4.1]{Zhang26} on $r\Disc$ to the
closed form $\Ric^{(1)}_\beta$ and the family $g_{\tau e_j}$. Since $V_j(0)=0$, we obtain
\begin{equation}
 \left.\partial_\tau\partial_{\bar\tau}\right|_{\tau=0}
 \frac1\pi\int_{r\Disc}\log\frac r{|z|}\,g_{\tau e_j}^*\Ric^{(1)}_\beta
 =\int_\T\Ric^{(1)}_\beta[V_j(r\zeta)]\,d\sigma(\zeta).
\end{equation}
By \Cref{v:l:factor} and the product rule, the density of
$g_{\tau e_j}^*\Ric^{(1)}_\beta$ is $C^2$ in $\tau$ with values in $L^2(\Disc)$.
The weights $\pi^{-1}\log(r/|z|)$, extended by zero off $r\Disc$, converge in
$L^2(\Disc)$ to $\pi^{-1}\log(1/|z|)$. Since $g$ and $V_j$ are continuous
on $\overline\Disc$, letting $r\uparrow1$ gives
\begin{equation}\label{eq2}
 \lim_{r\uparrow1}
 \left.\partial_\tau\partial_{\bar\tau}
 \left(\frac1\pi\int_{r\Disc}\log\frac r{|z|}\,
       g_{\tau e_j}^*\Ric^{(1)}_\beta\right)\right|_{\tau=0}
 =
 \left.\partial_\tau\partial_{\bar\tau}
 \mathcal P(g_{\tau e_j})\right|_{\tau=0}=\int_\T\Ric^{(1)}_\beta[V_j(\zeta)]\,d\sigma(\zeta).
\end{equation}
 Since $u(0)=\mathcal P(g)$, then \eqref{determiant}, \eqref{eq2} and  \Cref{v:l:factor} give
\begin{equation}\label{v:identities}
 \begin{aligned}
 &\log\det(B_{x_0}^*\eta(x_0)B_{x_0})=\mathcal P(g)-W_{x_0},\\
 &\mathcal L\mathcal P
 =\int_\T(\tr_\beta\Ric^{(1)}_\beta)(g(\zeta))\,d\sigma(\zeta)\ge -n .
 \end{aligned}
\end{equation}
In an $\eta(x_0)$-unitary basis, let $v$ be the coordinate column of $g'(0)$
and let $I_n$ be the $n\times n$ identity matrix.  Let $N_\Lambda(v)$ be the complex Hessian at $v$ of
$\Phi_\Lambda(w):=-(n+1)\log(1-|w|^2/\Lambda)$. Direct computation  gives
\begin{equation}\label{eq:central-ball-hessian}
 N_\Lambda(v)=\frac{n+1}{\Lambda-|v|^2}
 \left(I_n+\frac{vv^*}{\Lambda-|v|^2}\right),\qquad
 \det N_\Lambda(v)=\left(\frac{n+1}{\Lambda}\right)^n e^{q_\Lambda(g)}.
\end{equation}
By \eqref{v:eq:centre-variations}, the differential of the holomorphic map
$t\mapsto g_t'(0)$ at $0$ is $B_{x_0}$. 
By \eqref{v:identities}, \eqref{eq:central-ball-hessian},
\begin{equation}\label{lowerbounde}
 \begin{aligned}
 \mathcal L q_\Lambda
 &=\left.\sum_{j=1}^n\partial_{t_j}\partial_{\bar t_j}
       q_\Lambda(g_t)\right|_{t=0}=\sum_{j=1}^n B_j(0)^*N_\Lambda(v)B_j(0)
   =\tr\!\bigl(B_{x_0}^*N_\Lambda(v)B_{x_0}\bigr)\\
 &\ge n\bigl(\det(B_{x_0}^*N_\Lambda(v)B_{x_0})\bigr)^{1/n}=n\bigl(\det N_\Lambda(v)\,
          \det(B_{x_0}^*\eta(x_0)B_{x_0})\bigr)^{1/n}\\
 &=\frac{n(n+1)}{\Lambda}
   \exp\!\left(\frac{q_\Lambda(g)+\mathcal P(g)-W_{x_0}}{n}\right).
 \end{aligned}
\end{equation}
 Each $g_t$ has center $x_0$ and finite area, so
$F_\delta(g_t)\le F_\delta(g_0)=b(\delta)$. Thus, for each $j$,
$\tau\mapsto F_\delta(g_{\tau e_j})$ has a maximum at $\tau=0$, and then 
\[
 \mathcal L F_\delta
 =\left.\sum_{j=1}^n\partial_\tau\partial_{\bar\tau}
      F_\delta(g_{\tau e_j})\right|_{\tau=0}\le0 .
\]
Since $F_\delta=q_\Lambda+\mathcal P-\delta\mathcal A$, then   \eqref{eq:source-area-functional-hessian}, \eqref{v:identities}
and \eqref{lowerbounde} give
\begin{equation}\label{v:basic}
 \begin{aligned}
 \frac{n(n+1)}{\Lambda}
 \exp\!\left(\frac{q_\Lambda(g)+\mathcal P(g)-W_{x_0}}{n}\right)
 &\le \mathcal Lq_\Lambda\le -\mathcal L\mathcal P+\delta\mathcal L\mathcal A\le n+C_1\delta\bigl(1+A_h(g)\bigr).
 \end{aligned}
\end{equation}
Put $s_\delta:=\delta\mathcal A(g_\delta)$. Since
$q_\Lambda(g_\delta)+\mathcal P(g_\delta)=b(\delta)+s_\delta$
and $\delta A_h(g_\delta)\le C_0s_\delta$, we have
\begin{equation}\label{v:numerical}
 e^{(b(\delta)+s_\delta-W_{x_0})/n}
 \le\frac{\Lambda}{n+1}\left(1+\frac{C_1}{n}(\delta+C_0s_\delta)\right).
\end{equation}
Fix $\Lambda>D_{x_0}$. Multiplying \eqref{v:numerical} by
$e^{-s_\delta/n}$, we obtain
\begin{equation}
 e^{(b(\delta)-W_{x_0})/n}
 \le \frac{\Lambda}{n+1}
 \left(1+\frac{C_1}{2n}+\frac{C_1C_0}{e}\right).
\end{equation}
Thus $b(\delta)$ is bounded above independently of $\delta$.
For $0<\delta'<\delta$, every finite-area holomorphic disc $f$ with center $x_0$
satisfies
\begin{equation}
 F_{\delta'}(f)-F_\delta(f)=(\delta-\delta')\mathcal A(f)\ge0.
\end{equation}
Taking suprema gives $b(\delta')\ge b(\delta)$. Hence $b(\delta)$
increases to a finite limit as $\delta\downarrow0$.

Since $F_\delta(g_\delta)=b(\delta)$ and
$F_{\delta/2}(g_\delta)\le b(\delta/2)$,
\begin{equation}\label{eq:penalty-limit}
 0\le \frac{s_\delta}{2}
 =\frac{\delta}{2}\mathcal A(g_\delta)
 =F_{\delta/2}(g_\delta)-F_\delta(g_\delta)
 \le b(\delta/2)-b(\delta)\longrightarrow0,
\end{equation}
then $s_\delta$ tends to $0$ as $\delta\downarrow 0$.
Put $b_*:=\lim_{\delta\downarrow0}b(\delta)$. 
Letting $\delta\downarrow0$ in \eqref{v:numerical} gives
\begin{equation}
 e^{(b_*-W_{x_0})/n}\le\frac{\Lambda}{n+1}.
\end{equation}
Now fix a finite-area holomorphic disc $f$ with center $x_0$. Since $q_\Lambda(f)\ge0$,
\begin{equation}
 \mathcal P(f)
 \le F_\delta(f)+\delta\mathcal A(f)
 \le b(\delta)+\delta\mathcal A(f).
\end{equation}
As $\mathcal A(f)<\infty$, letting $\delta\downarrow0$ yields
\begin{equation}
 \mathcal P(f)\le b_*\le W_{x_0}+n\log\frac{\Lambda}{n+1}.
\end{equation}
This holds for every $\Lambda>D_{x_0}$. Letting $\Lambda\downarrow D_{x_0}$
and using $e^{W_{x_0}}=\beta^n(x_0)/\eta^n(x_0)$, we obtain
\eqref{eq:fixed-disc-volume}. For the constant disc $f\equiv x_0$,
we have $\mathcal P(f)=0$, so $D_{x_0}\le D$ gives \eqref{eq:vol}
at every point where $\eta(x_0)>0$. At all other points $\eta^n(x_0)=0$,
and \eqref{eq:vol} follows there as well.
\end{proof}

\begin{proof}[Proof of \Cref{v:t:main}]
Take $\eta=h$ and $D=D_h$ in \Cref{prop:relative-volume}.
The bound \eqref{eq:source-centre-bound} holds by the definition of
$D_h$, and \Cref{le:brody-derivative-bound} gives $0<D_h<\infty$.
The global comparison \eqref{eq:vol} is precisely \eqref{v:out}.
\end{proof}

\section{Ampleness of the canonical bundle}\label{sec:backend}\label{sec:ampleness}

\subsection{From pointwise volume comparison to canonical volume}

\begin{lemma}\label[lemma]{prop:nonnegative-density}
Let $M$ be a compact complex $n$-fold, and let $\nu\ge0$ be a continuous
$(n,n)$-form with $\nu\not\equiv0$. Suppose that every Hermitian metric
$\beta$ on $M$ satisfies
\begin{equation}
 \tr_\beta\Ric^{(1)}_\beta\ge-n\quad\Longrightarrow\quad\beta^n\ge\nu .
\end{equation}
For every Hermitian metric $h$ and every Gauduchon metric $\gamma$
$(\ddc\gamma^{n-1}=0)$, we have
\begin{equation}\label{eq:gauduchon-test}
 -\int_M\Ric^{(1)}_h\wedge\gamma^{n-1}
 \ge\int_M r_\gamma\gamma^n,\qquad
 r_\gamma:=\left(\frac{\nu}{\gamma^n}\right)^{1/n}.
\end{equation}
If  $\nu=h^n$ and $\gamma$ is conformal to $h$, equality
implies that $h$ is K\"ahler.
If $M$ is K\"ahler, then
\begin{equation}\label{eq:density-canonical-volume-bound}
 \vol(K_M)\ge(2\pi)^{-n}\int_M\nu>0 .
\end{equation}
\end{lemma}

\begin{proof}
\textbf{Step 1. The Gauduchon inequality.}
Fix a Gauduchon metric $\gamma$ and write $S_\beta:=\tr_\beta\Ric^{(1)}_\beta$.
For $0<2\varepsilon<\max_M r_\gamma$, we can find  a real-valued smooth function $q$ with $r_\gamma-2\varepsilon<q<r_\gamma$.
Set 
\[
c:=
\frac{\displaystyle\int_M(S_\gamma+nq)\,\gamma^n}
{\displaystyle\int_M\gamma^n},
\]
then $\int_{X}(S_\gamma+nq-c)\gamma^n=0$. Since $\gamma$ is a  Gauduchon metric, we can find 
a real-valued function $u\in C^\infty(M)$ such that
$\Delta_\gamma u=S_\gamma+nq-c$ \cite[Proposition~9, (22)]{GaudRic77},
\cite[proof of Theorem~3.1, (3.1)]{ACS17}. Since $\max_M q>0$, by adding  a constant to $u$,
we may assume $\max_M(qe^{-u/n})=1$. For $\beta:=e^{u/n}\gamma$,
\begin{equation}\label{eq:conformal-chern-scalar}
 \beta^n=e^u\gamma^n,\qquad
 S_\beta=e^{-u/n}(S_\gamma-\Delta_\gamma u)=e^{-u/n}(c-nq).
\end{equation}
If $c\ge0$, then $S_\beta\ge-n$ everywhere. Thus,  at a maximum $p$ of
$qe^{-u/n},$ we have $\beta^n(p)=q(p)^n\gamma^n(p)<\nu(p)$. This
contradicts the hypothesis. Thus $c<0$, and $\varepsilon\downarrow0$ gives
\begin{equation}\label{eq:gauduchon-scalar-density-bound}
 -\frac1n\int_M S_\gamma\gamma^n\ge\int_M r_\gamma\gamma^n .
\end{equation}
Set $\alpha:=-\Ric^{(1)}_h$, then integration by parts gives 
\begin{equation}
\int_M
\left(
-\operatorname{Ric}^{(1)}_\gamma-\alpha
\right)\wedge\gamma^{n-1}=
2\pi\int_M
\log\frac{\gamma^n}{h^n}\,
dd^c\gamma^{n-1}
=0,
\end{equation}
which implies 
\begin{equation}\label{equality11}
-\frac1n\int_MS_\gamma\,\gamma^n
=
-\int_M\operatorname{Ric}^{(1)}_\gamma\wedge\gamma^{n-1}\\
=
\int_M\alpha\wedge\gamma^{n-1}.
\end{equation}
Combining this equality with \eqref{eq:gauduchon-scalar-density-bound} proving
\eqref{eq:gauduchon-test}. 

For $n=1$, \eqref{equality11} gives $2\pi\deg K_M\ge\int_M\nu>0$, and Riemann--Roch theorem gives
$\vol(K_M)=\deg K_M$ \cite{Rie57,Roc65}, \cite[16.9 and~17.16]{Forster81}.

\textbf{Step 2. A positive current and the equality case.}
Assume $n\ge2$. Let $V$ be the Fr\'echet space of smooth real
$(n-1,n-1)$-forms. For each $\Omega\in V$ with $\Omega>0$, there is a unique
smooth Hermitian form $\gamma$ such that $\Omega=\gamma^{n-1}$ \cite[pp.~279--280, (4.8)]{Michelsohn82}. For such an $\Omega$, put
$J(\Omega):=\int_M r_\gamma\gamma^n$, then $J$ is continuous, concave and positively homogeneous.
 Pointwise, for $A\ge0$ with $A^n\ge\nu$,
\begin{equation}
 A\wedge\Omega=\frac{\tr_\gamma A}{n}\gamma^n
 \ge\left(\frac{A^n}{\gamma^n}\right)^{1/n}\gamma^n
 \ge r_\gamma\gamma^n,\qquad
 (r_\gamma\gamma)\wedge\Omega=r_\gamma\gamma^n.
\end{equation}
The choice $A:=r_\gamma\gamma$ attains the minimum. 

Write $\langle\alpha,\Omega\rangle:=\int_M\alpha\wedge\Omega$. Define 
\begin{equation}
 E_1:=\{(\Omega,t):\Omega>0,\ t<J(\Omega)\}\subset V\times\R,\qquad
 E_2:=\{(\Omega,\langle\alpha,\Omega\rangle):\ddc\Omega=0\}
 \quad\subset V\times\R,
\end{equation}
then $E_1$ is an open convex cone, and $E_2$ is a linear space.
By \eqref{eq:gauduchon-test}, $E_1\cap E_2=\emptyset$. Hahn--Banach separation gives a continuous functional
$(\Omega,t)\mapsto\langle R,\Omega\rangle-at$, which is strictly positive in $E_2$ and zero on $E_2$. Sending $t\to-\infty$ gives $a\ge0$.
For a Gauduchon power $\Omega=\gamma^{n-1}$ \cite[Theorem~1]{Gaud77},
the two sets contain points with the same first coordinate, so $a>0$.
Normalize $a=1$ and let $t\uparrow J(\Omega)$. Then
\begin{equation}
 \langle R,\Omega\rangle\ge J(\Omega)\ge0\quad(\Omega>0),\qquad
 \langle R-\alpha,\Omega\rangle=0\quad(\ddc\Omega=0).
\end{equation}
In particular $R\ge0$, and Lamari's duality
\cite[\S1.3(1), p.~267]{Lam99} gives a real distribution $\varphi$ with
\begin{equation}\label{eq:canonical-current-class}
 R=\alpha+2\pi\ddc\varphi,\qquad dR=0 .
\end{equation}
Write $R=R_{\rm ac}+R_{\rm s}$ for the Lebesgue decomposition of its
measure coefficients. The form $R_{\rm ac}$ is semipositive almost everywhere.
For any Hermitian metric $\gamma$ and smooth $\chi\ge0$, test with
$\Omega:=(\chi+\varepsilon)\gamma^{n-1}$. By pointwise homogeneity,
letting $\varepsilon\downarrow0$ gives
\begin{equation}\label{almost}
 \int_M\chi\,R\wedge\gamma^{n-1}\ge\int_M\chi\,r_\gamma\gamma^n,\qquad
 R_{\rm ac}\wedge\gamma^{n-1}\ge r_\gamma\gamma^n\quad\text{a.e.}
\end{equation}
Here the second inequality follows by taking the absolutely continuous
part of the measure inequality \cite[Ch.~III, Proposition~1.14]{DemCADG}.
Note that the last inequality of \eqref{almost} holds almost everywhere for each fixed Hermitian metric $\gamma$, but the exceptional null set may depend on the choice of $\gamma$.

Cover $M$ by finitely many coordinate balls $B'\Subset B$. On each $B'$
take a countable dense family of constant positive Hermitian matrices,
extending each to a global metric by interpolating with $h$ in $B\setminus B'$.
Discard the union of the null sets where the pointwise inequality fails for
one of these metrics. At every remaining point, continuity in the test matrix gives
$R_{\rm ac}\wedge H^{n-1}\ge(\nu/H^n)^{1/n}H^n$ for every $H>0$.
In particular, $H:=R_{\rm ac}+\varepsilon h$ satisfies
\begin{equation}
 H^n\ge R_{\rm ac}\wedge H^{n-1}\ge(\nu/H^n)^{1/n}H^n .
\end{equation}
Thus $H^n\ge\nu$, and $\varepsilon\downarrow0$ gives $R_{\rm ac}^n\ge\nu$
almost everywhere.

Suppose equality holds in \eqref{eq:gauduchon-test}, with $\nu=h^n$ and
$\gamma$ conformal to $h$, so $r_\gamma\gamma=h$. Since $\ddc\gamma^{n-1}=0$,
the positive measure $R\wedge\gamma^{n-1}-r_\gamma\gamma^n$ has mass zero.
Thus $R_{\rm s}\wedge\gamma^{n-1}=0$ and positivity gives $R_{\rm s}=0$.
For the eigenvalues $a_i$ of $R_{\rm ac}$ relative to $\gamma$ we now have
\begin{equation}
 r_\gamma=\tfrac1n\sum_i a_i\ge(\textstyle\prod_i a_i)^{1/n}\ge r_\gamma
 \quad\text{a.e.},
\end{equation}
where the last inequality is $R_{\rm ac}^n\ge\nu$. Equality forces every
$a_i=r_\gamma$, hence $R=h$ and $dh=dR=0$.
\textbf{Step 3. The canonical volume in the K\"ahler case.}
Suppose $M$ is K\"ahler. By \eqref{eq:canonical-current-class},
$T:=R/(2\pi)$ is a closed positive representative of $c_1(K_M)$,
which means that $K_M$ is pseudoeffective. Boucksom's estimate \cite[Proposition~3.1]{Bou02} gives
\begin{equation}
 \vol(K_M)\ge\int_M T_{\rm ac}^n
 =(2\pi)^{-n}\int_M R_{\rm ac}^n
 \ge(2\pi)^{-n}\int_M\nu .
\end{equation}
For the gluing in \cite[proof of Theorem~2.4, pp.~1050--1051]{Bou02}, write
$T=\theta+\ddc\varphi$ with $\theta=\alpha/(2\pi)$ and $\varphi$ quasi-plurisubharmonic
(locally a plurisubharmonic function plus a smooth function). Write $\nu_x(T)$
for the Lelong number, given by $\liminf_{z\to0}u(z)/\log|z|^2$ for a local plurisubharmonic
potential $u$ with $T=\ddc u$ in coordinates centered at $x$. Let $d_k\ge\varphi$ be the Bergman weights,
$u_j\downarrow\varphi$ the smooth regularizations, and $b_k:=kd_k/(k+1)-k$.
Choose $\{d_k=-\infty\}\subset W_k\Subset U_k\Subset\{d_k<-(k+1)^2\}$.
On $\overline U_k\setminus W_k$, $b_k$ is continuous and $b_k>d_k+1\ge\varphi+1$.
Thus $u_j<b_k-1/2$ there for all large $j$. Since
$\sup_{M\setminus W_k}\nu_x(T)\to0$ \cite[Proposition~3.7(ii)]{DemReg92}, choose such
$j_k\ge k$ with $\theta+\ddc u_{j_k}\ge-\varepsilon_k\omega$ off $W_k$,
where $\omega$ is K\"ahler and $\varepsilon_k\to0$, as in Boucksom's proof.
The weight $\psi_k=b_k$ on $U_k$ and $\psi_k=\max_{\rm reg}(b_k,u_{j_k})$
off $\overline W_k$ is well defined if the regularized maximum agrees with
the larger entry whenever the gap exceeds $1/2$.
By \cite[Proposition~3.7(iii)]{DemReg92}, $\theta+\ddc d_k\ge-\eta_k\omega$
with $\eta_k\to0$. Thus $\theta+\ddc b_k=\tfrac{k}{k+1}(\theta+\ddc d_k)
+\tfrac{1}{k+1}\theta\ge-o(1)\omega$ on $M$. Off $\overline W_k$, the
regularized maximum preserves the common lower curvature bound of its two
entries, so $\theta+\ddc\psi_k\ge-o(1)\omega$ on $M$. At each $x$ with
$\varphi(x)>-\infty$, $d_k(x)\to\varphi(x)$, so eventually
$x\notin\overline U_k$ and $b_k(x)<u_{j_k}(x)-1/2$. Thus $\psi_k=u_{j_k}$
near $x$, giving $(\theta+\ddc\psi_k)_{\rm ac}(x)\ge T_{\rm ac}(x)-t\omega(x)$
eventually almost everywhere for each $t>0$
\cite[Lemma~8.6 and Theorem~9.1]{DemReg82}. This is the lower bound used
with Fatou's lemma on the locus where $T_{\rm ac}$ is positive definite.
\end{proof}

\subsection{Canonical ampleness and projectivity}\label{sec:canonical-ampleness}

\begin{restatement}{\Cref{thm:main}}
Let $X$ be a compact hyperbolic complex manifold. Then $K_X$ is ample
and $X$ is projective.
\end{restatement}

\begin{proof}
Fix a Hermitian metric $h$ on $X$ and put $n:=\dim X$. For $x_0\in X$, set
\begin{equation}
 D_h(x_0):=\sup_f |f'(0)|_h^2,
\end{equation}
where $f$ ranges over holomorphic discs in $X$ with center $x_0$.
Small coordinate discs and \Cref{le:brody-derivative-bound} give
$0<D_h(x_0)\le D_h<\infty$. The function $x_0\mapsto D_h(x_0)$ is upper
semicontinuous. Indeed, for $x_j\to x_0$ take a subsequence realizing the
upper limit and discs $f_j$ centered at $x_j$ with $D_h(x_j)\le |f_j'(0)|_h^2+1/j$.
By \eqref{g:eq:brody-speed} and compactness of $X$, a further subsequence
converges locally smoothly to a holomorphic disc $f$ with center $x_0$, giving
$\limsup_jD_h(x_j)=|f'(0)|_h^2\le D_h(x_0)$. Thus
$\nu(x_0):=((n+1)/D_h(x_0))^n h^n(x_0)$ is measurable. Every Hermitian metric
$\beta$ with $\tr_\beta\Ric^{(1)}_\beta\ge-n$ satisfies $\beta^n\ge\nu$
by \eqref{eq:fixed-disc-volume} applied to the constant disc. We prove that
\begin{equation}\label{eq:canonical-volume-current}
 \vol(K_X)\ge\int_X\left(\frac{n+1}{2\pi D_h(x_0)}\right)^n h^n
 \ge\left(\frac{n+1}{2\pi D_h}\right)^n\int_X h^n>0 .
\end{equation}

\textbf{Step 1. A bounded potential in the canonical class.}
Put $\alpha:=-\Ric^{(1)}_h$ and let $\mathcal P_h$ denote the Green
integral defined before \Cref{v:l:radial}, with $\beta=h$.
With $a:=\max\{1,-\tfrac1n\min_X\tr_h\Ric^{(1)}_h\}$ we have
$\tr_{ah}\Ric^{(1)}_{ah}=a^{-1}\tr_h\Ric^{(1)}_h\ge-n$. Taking
$\eta:=h$ and $\beta:=ah$ in \eqref{eq:fixed-disc-volume} gives, for
every finite-area disc $f$ centered at $x_0$,
\begin{equation}
 \mathcal P_h(f)\le n\log\frac{aD_h(x_0)}{n+1}\le G,\qquad
 G:=n\log\frac{aD_h}{n+1}.
\end{equation}
Let $b(x_0)$ be the infimum of $-\mathcal P_h(f)$ over closed holomorphic
discs centered at $x_0$. The preceding estimate and the constant disc give
$-G\le b\le0$. For the smooth closed form $\alpha=-\Ric^{(1)}_h$, the disc
envelope in \cite[Theorem~1.2 and \S7]{LSWZ26} is $b/2$, and the operator
$dd^c$ there is $4\pi\ddc$ in our notation. The theorem therefore gives
$\alpha+2\pi\ddc b\ge0$. With $\pi\colon(\det T_X)^\times\to X$ the bundle projection, the Poisson
envelope $E$ of the log squared norm
\cite[\S7, Theorem]{Rosay03} satisfies $E(v)=\log|v|_{\det h}^2+b(\pi(v))$.
Indeed, each closed disc $f$ centered at $\pi(v)$ admits a nonvanishing
lift with center $v$, since the pullback line bundle is holomorphically
trivial on a slightly larger disc \cite[Theorem~30.4]{Forster81}.
By \eqref{eq:disc-jensen}, every such lift has boundary mean
$\log|v|_{\det h}^2-\mathcal P_h(f)$ for its log squared norm. Since $b$ is bounded, the positive current
has zero Lelong numbers, and $[\alpha]$ is nef by \cite[Corollary~6.4]{DemReg92}.

\textbf{Step 2. Smooth metrics with uniformly bounded potentials.}
For each $s>0$, nefness gives $\psi\in C^\infty(X,\R)$ with
$\alpha+2\pi\ddc\psi\ge-h/(2s)$. Thus
\begin{equation}
 h+s\alpha+2\pi\ddc(s\psi)\ge h/2 .
\end{equation}
The existence interval in \cite[Theorem~1.1]{ShermanWeinkove20}
is therefore $[0,\infty)$. Let $\omega(s)>0$ solve
$\omega(s)=h-s\Ric^{(1)}_{\omega(s)}$ and put
$\beta_\varepsilon:=\varepsilon\omega(1/\varepsilon)$,
$u_\varepsilon:=\log(\beta_\varepsilon^n/h^n)$ for $0<\varepsilon\le1$.
Since constant rescaling preserves the Chern--Ricci form,
\begin{equation}
 \beta_\varepsilon
 =\varepsilon h-\Ric^{(1)}_{\beta_\varepsilon}
   =\varepsilon h+\alpha+2\pi\ddc u_\varepsilon,\qquad
 \tr_{\beta_\varepsilon}\Ric^{(1)}_{\beta_\varepsilon}
 =\varepsilon\tr_{\beta_\varepsilon}h-n\ge-n .
\end{equation}
The trace bound gives $\beta_\varepsilon^n\ge\nu$ and thus
$e^{u_\varepsilon}\ge((n+1)/D_h)^n$. Put
$C:=1+\max_{|v|_h=1}|\alpha[v]|$, so $\varepsilon h+\alpha\le Ch$
for $0<\varepsilon\le1$. At a maximum of $u_\varepsilon$,
$\ddc u_\varepsilon\le0$ gives $\beta_\varepsilon\le Ch$.
We obtain
\begin{equation}
 n\log((n+1)/D_h)\le u_\varepsilon\le n\log C\qquad\text{on }X .
\end{equation}

\textbf{Step 3. The volume estimate.}
On coordinate balls $B'\Subset B$ relatively compact in a chart, choose $A$
with $\alpha+\varepsilon h,h\le A\omega_{\rm e}$ on $B$ for $0<\varepsilon\le1$,
where $\omega_{\rm e}=2\pi\ddc|z|^2$. Then
\begin{equation}
 U_\varepsilon:=u_\varepsilon+A|z|^2,\qquad
 0<\beta_\varepsilon\le2\pi\ddc U_\varepsilon,\qquad
 \sup_\varepsilon\|U_\varepsilon\|_{L^\infty(B)}<\infty .
\end{equation}
The $U_\varepsilon$ are plurisubharmonic. The Chern--Levine--Nirenberg estimate
\cite{CLN69}, \cite[Ch.~III, (3.3)]{DemCADG} gives
\begin{equation}
 \sup_\varepsilon\int_{B'} h\wedge\beta_\varepsilon^{n-1}
 \le A(2\pi)^{n-1}\sup_\varepsilon
       \int_{B'}\omega_{\rm e}\wedge(\ddc U_\varepsilon)^{n-1}<\infty .
\end{equation}
For $n=1$ the local integral of $h$ is fixed. A finite cover yields
$\int_X h\wedge\beta_\varepsilon^{n-1}\le C'$.
Set $\theta_\varepsilon:=\beta_\varepsilon-\varepsilon h
=-\Ric^{(1)}_{\beta_\varepsilon}$. The metric induced by
$\beta_\varepsilon$ on $K_X$ has Chern form $\theta_\varepsilon/(2\pi)$.
Let $G_\varepsilon$ consist of the points where $\theta_\varepsilon$
has no zero eigenvalues and at most one negative eigenvalue relative to $h$.
Applied to $\beta_\varepsilon$ and $\varepsilon h$, Demailly's pointwise inequality
\cite[Lemma~9.46, $q=1$]{Dem20} gives
\begin{equation}
 \mathbf1_{G_\varepsilon}\theta_\varepsilon^n
 \ge\beta_\varepsilon^n-n\varepsilon h\wedge\beta_\varepsilon^{n-1}.
\end{equation}
For fixed $\varepsilon$, the strong holomorphic Morse inequality in
degree one \cite[Theorem~0.1(b)]{DemMorse85}, \cite[\S5]{Berman04} gives
\begin{equation}
 (2\pi)^n\vol(K_X)\ge\int_{G_\varepsilon}\theta_\varepsilon^n
 \ge\int_X\beta_\varepsilon^n
   -n\varepsilon\int_X h\wedge\beta_\varepsilon^{n-1}
 \ge\int_X\nu-nC'\varepsilon .
\end{equation}
Letting $\varepsilon\downarrow0$ proves \eqref{eq:canonical-volume-current}.
In particular $K_X$ is big. The section growth bound
$h^0(X,mK_X)=O(m^{a(X)})$, where $a(X)$ is algebraic dimension of $X$. Thus $a(X)=n$ and $X$ is a Moishezon manifold \cite[Theorem~5.1]{DemMorse85}.
 Since $X$ contains no rational curves,  $X$ is
projective by \cite[Corollary~1.4.6]{BCHM10} (see \cite[Theorem~4.4]{HM10} for the smooth Moishezon case). Since $X$ is of general type, projective and contains no ration curves, 
$K_X$ is ample by \cite[Lemma~2.1]{DT}.
\end{proof}

\subsection{Constant negative Chern holomorphic sectional curvature}

\begin{restatement}{\Cref{cor:intro-constant-hsc}}
Let $h$ be a Hermitian metric of constant Chern holomorphic sectional
curvature $c<0$ on a compact complex manifold $X$. Then $h$ is
K\"ahler and $(X,h)$ is a complex space form, that is, a complete connected K\"ahler manifold of constant holomorphic sectional curvature.
\end{restatement}

\phantomsection\label{proof:constant-negative-hsc}
\begin{proof}
Put $n:=\dim X$. When $n=1$, $dh=0$ by degree and $h$ is complete by
compactness, proving the claim. Assume $n\ge2$ and rescale $h$ so that
$c=-2/(n+1)$. The disc metric of this curvature gives, by Kobayashi's
Schwarz lemma \cite[Theorems~3.8 and~9.2]{Kob67}, hyperbolicity and
\begin{equation}\label{eq:hermitian-disc-schwarz}
 q^*h\le(n+1)\frac{\sqrt{-1}\,dz\wedge d\bar z}{(1-|z|^2)^2}
 \qquad(q\colon\Disc\to X\text{ holomorphic}).
\end{equation}
Evaluating at $0$ gives $D_h\le n+1$. Thus \Cref{v:t:main} gives
$\beta^n\ge((n+1)/D_h)^nh^n\ge h^n$ for every Hermitian $\beta$ with
$\tr_\beta\Ric^{(1)}_\beta\ge-n$. Apply \Cref{prop:nonnegative-density}
with $\nu:=h^n$ and a Gauduchon metric $\gamma:=v^{1/(n-1)}h$, where
$v>0$ is smooth \cite[Theorem~1]{Gaud77}. Then
\begin{equation}
 \gamma^{n-1}=v h^{n-1},\qquad
 r_\gamma=(\nu/\gamma^n)^{1/n}=v^{-1/(n-1)},\qquad
 r_\gamma\gamma=h,\quad r_\gamma\gamma^n=vh^n.
\end{equation}
Put $s:=\tr_h\Ric^{(1)}_h$ and write $\partial_\gamma^*$ for the formal
adjoint of $\partial$ with respect to $\gamma$. \cite[(3.1), (3.5), (3.8)]{Yang16} give
\begin{equation}\label{eq:rigidity-weighted-scalars}
 -\int_X\Ric^{(1)}_h\wedge\gamma^{n-1}
 =-\frac1n\int_X vs\,h^n
 =-\frac{(n+1)c}{2}\int_X vh^n
   -\frac1{2n}\int_X|\partial_\gamma^*\gamma|_\gamma^2\,\gamma^n
 \le\int_X r_\gamma\gamma^n ,
\end{equation}
where the last step uses $c=-2/(n+1)$. By \eqref{eq:gauduchon-test}
equality holds, and its conformal equality case gives $dh=0$. The original metric is
therefore K\"ahler as well. It has constant HSC and is complete by
compactness, hence is a complex space form.
\end{proof}

\section{General type of subvarieties}\label{sec:general-type}

\subsection{Canonical volume and generically finite maps}\label{sec:plurigenera}

A modification is a proper holomorphic map that restricts to a biholomorphism
between the complements of nowhere dense analytic subsets of the source and target.
A K\"ahler modification is a modification from a K\"ahler manifold.
\begin{corollary}\label[corollary]{cor:smooth-semipositive-estimates}
Let $P$ be a compact complex $n$-fold in Fujiki's class $\mathcal C$,
$n\ge1$, and let $\eta$ be a smooth semipositive $(1,1)$-form with
$\int_P\eta^n>0$. If $|f'(0)|_\eta^2\le D$ for every holomorphic
disc $f\colon\Disc\to P$, where $D>0$, then
\begin{equation}\label{eq:smooth-semipositive-volume}
 \vol(K_P)\ge\left(\frac{n+1}{2\pi D}\right)^n\int_P\eta^n .
\end{equation}
\end{corollary}
\begin{proof}
Take a K\"ahler modification $\rho\colon Q\to P$
\cite[Chapter~IV, \S3.2, Theorem~5]{Var89}. Composition with $\rho$ gives
$|f'(0)|_{\rho^*\eta}^2\le D$ for every disc in $Q$.
Thus \Cref{prop:relative-volume,prop:nonnegative-density}, with
$\nu=((n+1)/D)^n(\rho^*\eta)^n$, give \eqref{eq:smooth-semipositive-volume}
on $Q$. Since $\deg\rho=1$, integration of top forms gives
$\int_Q(\rho^*\eta)^n=\int_P\eta^n$, even if $d\eta\ne0$.
The plurigenera are bimeromorphic invariants \cite[Corollary~1.3]{Ueno73},
so $\vol(K_Q)=\vol(K_P)$.
\end{proof}

A resolution of an irreducible space $Z$ is a modification
$\mu\colon P\to Z_{\rm red}$ with $P$ smooth. Such resolutions exist in the
complex analytic category \cite[Theorem~13.4]{BM97}. Since $P$ is normal and
$\mu$ is dominant, it factors uniquely through the normalization
$\nu\colon Z^\nu\to Z_{\rm red}$, the finite bimeromorphic map with normal source
\cite[Ch.~VIII, \S4, no.~3]{GR84}.
The lift $P\to Z^\nu$ is proper since the preimage of any compact $K\subset Z^\nu$
is closed in $\mu^{-1}(\nu(K))$. It is biholomorphic over the common dense
isomorphism locus in $(Z_{\rm red})_{\rm reg}$, hence is a resolution of $Z^\nu$.
Bimeromorphic invariance of the plurigenera gives the same canonical volume
on every resolution. Resolutions of spaces in class
$\mathcal C$ belong to $\mathcal C$ \cite[\S4.3, Lemma~4.6(3)]{Fujiki78}.

\begin{proof}[Proof of \Cref{thm:proj}]
The reduced image $Y=g(Z)$ is closed analytic by
Remmert's projection theorem applied to the graph of $g$
\cite[\S7, Satz~20]{Remmert56} (see also \cite[Ch.~II, Theorem~8.8]{DemCADG}
for the reduced-space formulation). Since $Z$ is irreducible, so is $Y$. It inherits
hyperbolicity from $X$ and lies in $\mathcal C$ as a holomorphic image
of $Z$ \cite[\S4.3, Lemma~4.6(1)]{Fujiki78}. The fibre dimension
formula gives $\dim Y=\dim Z=:n$, since $g$ has finite general fibre
\cite[Ch.~II, (8.3) and Corollary~8.6]{DemCADG}.

Fix a Hermitian metric $h$ on $Y$ from local embeddings and a coordinate ball
$U\Subset Y_{\rm reg}$. Hyperbolicity gives
$0<C:=\sup_{q(0)\in U}|q'(0)|_h<\infty$, where $q\colon\Disc\to Y$
ranges over holomorphic discs \cite[Definition~1.1 and Theorem~1.2]{ThaiDuc00}.
Choose a nonzero smooth cutoff $0\le\chi\le1$ with compact support
in $U$ and extend
$\tau=C^{-2}\chi h$ by zero. For every generically finite surjection
$F\colon M\to Y$ from a compact manifold in $\mathcal C$, the form
$F^*\tau$ is smooth and semipositive, and
$|f'(0)|_{F^*\tau}^2\le1$ for every disc in $M$ by composition.
The restriction $F^{-1}(U)\to U$ is proper, with smooth source and target
and generic degree $\deg F$, the number of points in a general fibre. Applying the degree formula on each connected
component of the source gives $F_*[F^{-1}(U)]=(\deg F)[U]$
\cite[Ch.~III, (9.1) and Example~9.2]{DemCADG}.
Testing against $\tau^n$, which is compactly supported in $U$, gives
$\int_M(F^*\tau)^n=\deg(F)\int_U\tau^n$.
Thus \eqref{eq:smooth-semipositive-volume} gives
\begin{equation}
 \vol(K_M)\ge\deg(F)c_Y,\qquad
 c_Y:=\left(\frac{n+1}{2\pi}\right)^n\int_U\tau^n>0 .
\end{equation}
Take resolutions $\mu\colon P\to Z$ and $\sigma\colon Q\to Y$, both
with source in $\mathcal C$.
Since $\mu$ and $\sigma$ are modifications, $\deg(g\circ\mu)=\deg g$
and $\deg\sigma=1$. Taking $F=g\circ\mu$ and $F=\sigma$ gives
$\vol(K_P)\ge\deg(g)c_Y$ and $\vol(K_Q)\ge c_Y$.
Hence $Z$ and $Y$ are of general type.
\end{proof}

\begin{proof}[Proof of \Cref{cor:analytic-subspaces}]
By \Cref{thm:main}, $X$ is projective. Its closed analytic subspace
$Z_{\rm red}$ therefore belongs to class $\mathcal C$
\cite[\S4.3, Lemma~4.6(2)]{Fujiki78}. Applying \Cref{thm:proj} to the
degree-one inclusion $Z_{\rm red}\hookrightarrow X$ gives general type
for $Z_{\rm red}$ and hence for $Z$ by definition.
\end{proof}

\subsection{Canonical volume and projective degree}\label{sec:projective-degree}

Fix a projective embedding $\iota\colon X\hookrightarrow\mathbb P^N$
of a positive-dimensional hyperbolic complex space. Write
$L:=\iota^*\mathcal O_{\mathbb P^N}(1)$ and normalize
$[\omega_{\rm FS}]=2\pi c_1(\mathcal O_{\mathbb P^N}(1))$. Then
\begin{equation}\label{eq:embedded-disc-bound}
 0<D_X:=\sup_{q\colon\Disc\to X\ {\rm holomorphic}}
             |(\iota\circ q)'(0)|_{\rm FS}^2<\infty .
\end{equation}
Indeed, every disc factors through $X_{\rm red}$. The local derivative
bound of \cite[Definition~1.1 and Theorem~1.2]{ThaiDuc00} for the
induced Fubini--Study metric gives a global bound by a finite cover.
A nonconstant coordinate disc with nonzero derivative at a smooth point
of a positive-dimensional component of $X_{\rm red}$ gives $D_X>0$.
For an irreducible subvariety $Y$, write
$\deg_L Y=c_1(L)^{\dim Y}\cdot[Y]$.

\begin{corollary}\label[corollary]{cor:volume-degree}
Under the hypotheses of \Cref{thm:proj}, suppose $X$ has the fixed
projective embedding. Put $n:=\dim Z$ and $e:=\deg(g\colon Z\to Y)$.
For every resolution $\mu\colon P\to Z$, the map $p:=g\circ\mu$ satisfies
\begin{equation}\label{eq:volume-degree}
 \vol(K_P)\ge\left(\frac{n+1}{D_X}\right)^n(p^*L)^n,
 \qquad (p^*L)^n=\int_P c_1(p^*L)^n=e\,\deg_L Y .
\end{equation}
\end{corollary}
\begin{proof}
By Chow's theorem $\iota(Y)$ is algebraic \cite[Theorem~V]{Chow49}.
For $k_p:=(\iota\circ p)^*\omega_{\rm FS}$ the degree formula
$(\iota\circ p)_*[P]=e[\iota(Y)]$ \cite[Ch.~III, Example~9.2]{DemCADG} gives
\begin{equation}
 (2\pi)^{-n}\int_P k_p^n=(p^*L)^n=e\,\deg_L Y>0 .
\end{equation}
For every disc $f$ in $P$, \eqref{eq:embedded-disc-bound} applied to $p\circ f$
gives $|f'(0)|_{k_p}^2\le D_X$. Now apply \eqref{eq:smooth-semipositive-volume}.
\end{proof}

\subsection{Mixed intersections and the relative canonical model}\label{sec:relative-intersection}

\begin{lemma}\label[lemma]{lem:relative-intersection}
Let $p\colon P\to X$ be a generically finite map from a smooth
projective $n$-fold to a hyperbolic complex space with the fixed
projective embedding. Let $Y$ be the reduced image and
$e:=\deg(p\colon P\to Y)$. For any nef real $(1,1)$-classes
$\theta_1,\ldots,\theta_{n-1}$ on $P$,
\begin{equation}\label{eq:relative-intersection}
 c_1(K_P)\cdot\prod_{j=1}^{n-1}\theta_j
 \ge \frac{n+1}{D_X}(e\,\deg_L Y)^{1/n}
            \prod_{j=1}^{n-1}(\theta_j^n)^{1/n}.
\end{equation}
\end{lemma}
Here $c_1(K_P)\cdot\prod_j\theta_j$ is the integral of the wedge product
of smooth closed representatives of these classes on $P$.
\begin{proof}
By \eqref{eq:volume-degree}, $\vol(K_P)\ge((n+1)/D_X)^n e\,\deg_L Y>0$.
For $n=1$ the products are empty and $\deg K_P=\vol(K_P)$, proving
the claim. For $n\ge2$, apply Fujita approximation \cite[p.~1, Theorem]{Fujita94} to the big line
bundle $K_P$, in the ample form of \cite[Theorem~1.1]{Bou02}. For each $\varepsilon>0$ it gives a
smooth projective modification $\pi\colon\widetilde P\to P$ and
$\pi^*K_P=A+E$ as rational divisor classes, with $A$ ample, $E$
effective and $|A^n-\vol(K_P)|<\varepsilon$.
The classes $\pi^*\theta_j$ are nef, so their product has nonnegative
intersection with $E$ by approximation with K\"ahler classes. By the
projection formula and the mixed intersection inequality for nef real
$(1,1)$-classes \cite[Proposition~5.2(b), author's version]{Dem93},
\begin{equation}\label{eq:relative-intersection-approximation}
 c_1(K_P)\cdot\prod_j\theta_j
 \ge c_1(A)\cdot\prod_j\pi^*\theta_j
 \ge (A^n)^{1/n}\prod_j(\theta_j^n)^{1/n}.
\end{equation}
Letting $\varepsilon\downarrow0$ and using \eqref{eq:volume-degree}
gives \eqref{eq:relative-intersection}.
\end{proof}

For a projective morphism $g_Y\colon Z\to Y$ with $Z$ smooth, the relative
canonical algebra and the canonical ring of $Z$ are
\[
 \mathcal R(Z/Y):=\bigoplus_{m\ge0}(g_Y)_*\mathcal O_Z(mK_Z),\qquad
 R(Z):=\bigoplus_{m\ge0}H^0(Z,mK_Z).
\]
In the setting below, the \emph{relative canonical model} of $Z$ over $Y$ is
$\operatorname{Proj}_Y\mathcal R(Z/Y)$, and its \emph{absolute canonical model}
is $\operatorname{Proj}R(Z)$. A real Cartier divisor class on a projective variety
is the numerical class of a real linear combination of Cartier divisors. It is nef
if its intersection with every irreducible curve is nonnegative.

\begin{corollary}\label[corollary]{cor:model-ample}
In \Cref{cor:volume-degree}, suppose $Z$ is a projective manifold
and write $g_Y\colon Z\to Y$ for the map onto its reduced image.
The relative canonical model $V$ exists and is isomorphic to the absolute
canonical model of $Z$. Its canonical divisor $K_V$ is ample. For $m>0$
with $mK_V$ Cartier, write $c_1(K_V):=m^{-1}c_1(\mathcal O_V(mK_V))$. For any nef real
Cartier divisor classes $a_1,\ldots,a_{n-1}$ on $V$,
\begin{equation}\label{eq:relative-model-intersection}
 c_1(K_V)\cdot a_1\cdots a_{n-1}\ge
 \frac{n+1}{D_X}(e\,\deg_L Y)^{1/n}\prod_{j=1}^{n-1}(a_j^n)^{1/n}.
\end{equation}
\end{corollary}
\begin{proof}
Chow's theorem for the graph makes $g_Y$ an algebraic morphism
\cite[Theorems~V--VII]{Chow49}. Every morphism from an abelian variety
to $Y$ is constant by hyperbolicity \cite[Remark~2.6]{JK18}.
By \Cref{cor:volume-degree} and bimeromorphic invariance of the
plurigenera, every projective variety birational to $Z$ has Kodaira
dimension $n$ on a resolution. For a line bundle $L$ on a compact complex manifold $M$, its Iitaka dimension
$\kappa(M,L)$ is the largest dimension of the images of the meromorphic maps
defined by the nonempty systems $|L^{\otimes m}|$, $m\ge1$, and is $-\infty$
if all these systems are empty. We write $\kappa(M):=\kappa(M,K_M)$. In particular every projective model birational to $Z$ with canonical
singularities has Kodaira dimension $n>0$. This verifies the conclusion of
Hu--Meng--Zhang's Conjecture~1.3(2) on every model to which it applies, as required by
\cite[Proposition~2.2(5) and Theorem~3.1]{HMZ}. The theorem gives a model
$g_V\colon V\to Y$ birational to $Z$ over $Y$ with canonical singularities.
This means that $V$ is normal, $K_V$ is $\mathbb Q$-Cartier and
$K_P-\sigma_V^*K_V\ge0$ on every resolution $\sigma_V\colon P\to V$
with compatible canonical divisors. Cases~(2) and~(3) of that theorem require
$\kappa(V)=0$ and $\kappa(Z)\le n-3$, respectively. Hence case~(1) gives
$K_V$ ample, that is, $mK_V$ is an ample Cartier divisor for some $m>0$.

Take a common projective resolution $\sigma_T\colon P\to T$ over $Y$,
for $T=Z,V$. Smoothness of $Z$ and canonicity of $V$ give
$K_P=\sigma_T^*K_T+E_T$ with $E_T\ge0$ and
$\operatorname{codim}_T\sigma_T(\operatorname{Supp}E_T)\ge2$.
By normality and the projection formula,
$(\sigma_T)_*\mathcal O_P(mK_P)=\mathcal O_T(mK_T)$ for sufficiently divisible
$m$, compatibly with multiplication. Pushing forward to $Y$ or taking global
sections therefore identifies the respective graded algebras in these degrees.
Since $K_V$ is ample, both Proj constructions recover $V$
\cite[Theorem~3.52(1) and Corollary~3.53]{KM}.
Apply \eqref{eq:relative-intersection} on $P$ to $\sigma_V^*a_j$. The map
$P\to Y$ still has degree $e$. As $(\sigma_V)_*K_P=K_V$, the projection
formula gives \eqref{eq:relative-model-intersection}.
\end{proof}

\subsection{Log canonical and semi-log canonical pairs}\label{sec:slc-pairs}

For a normal pair $(T,B)$ with $B\ge0$ and $K_T+B$ real Cartier, let
$a(E,T,B):=\operatorname{ord}_E(K_P-\mu^*(K_T+B))$ for every prime divisor
$E$ on a smooth proper birational model $\mu\colon P\to T$, with compatible canonical divisors.
The pair is \emph{log canonical} if $a(E,T,B)\ge-1$ for all such $E$.
The reduced images $\mu(E)$ with $a(E,T,B)=-1$ are its log canonical centres.
A normal variety $T$ has log canonical singularities if $(T,0)$ is log canonical.

Let $W$ be reduced and equidimensional, with only ordinary double crossings
at its nonsmooth points in codimension one. Assume $W$ satisfies $S_2$, that is,
$\operatorname{depth}\mathcal O_{W,w}\ge\min\{2,\dim\mathcal O_{W,w}\}$ for every $w\in W$.
For the normalization $\sigma\colon W^\nu\to W$, let $C$ be the reduced
conductor divisor on $W^\nu$ over the codimension-one nonnormal locus of $W$.
A real boundary $\Delta=\sum_i d_i\Delta_i$ has distinct prime components and $0\le d_i\le1$.
Assume that $\Delta$ shares no component with $\sigma(C)$. The pair $(W,\Delta)$ is
\emph{semi-log canonical} (slc) if $K_W+\Delta$ is real Cartier and
$(W^\nu,\Delta^\nu+C)$ is log canonical, where $\Delta^\nu$ is the strict transform.
Its slc strata are the irreducible components of $W$ and the reduced $\sigma$-images
of the log canonical centres of $(W^\nu,\Delta^\nu+C)$
\cite[p.~201 and Definitions~2.1--2.5]{Fuj14}.
A real Cartier divisor is ample if its numerical class is a positive real
linear combination of ample Cartier classes.
\begin{corollary}\label[corollary]{cor:slc-ample}
Let $(W,\Delta)$ be a positive-dimensional projective slc pair with
real boundary and $W$ hyperbolic. Then $D:=K_W+\Delta$ is ample.
Fix a projective embedding of $W$ with $L$ and $D_W$ as in
\eqref{eq:embedded-disc-bound}. For an slc stratum $S$ of dimension
$d>0$, put $B:=\nu^*(D|_S)$ on its normalization $\nu\colon N\to S$.
Then $(D|_S)^d\ge((d+1)/D_W)^d\deg_L S$, and for any nef real Cartier
divisor classes $a_1,\ldots,a_{d-1}$ on $N$,
\begin{equation}\label{eq:slc-stratum-mixed}
 B\cdot a_1\cdots a_{d-1}\ge
 \frac{d+1}{D_W}(\deg_L S)^{1/d}\prod_{j=1}^{d-1}(a_j^d)^{1/d}.
\end{equation}
\end{corollary}
\begin{proof}
Hyperbolicity excludes nonconstant morphisms $\mathbb A^1_\C\to W$
and hence rational curves. Thus $D$ is nef by the slc cone theorem
\cite[Theorem~1.19(1)]{Fuj14}, and $B$ is nef on $N$. The finite map
$N\to S$ makes $N$ projective. Fix an ample Cartier divisor $H$ on $N$.
For each $t>0$, subadjunction \cite[Corollary~1.4]{FujSub22} gives an
effective real divisor $\Theta_t$ with
$K_N+\Theta_t\sim_{\mathbb R}B+tH=:A_t$.

Take a projective resolution $\mu\colon P\to N$. The composite
$P\to N\to S\hookrightarrow W$ has degree one onto $S$, so
\eqref{eq:relative-intersection} on $P$ has constant
$c_S:=(d+1)(\deg_L S)^{1/d}/D_W$. For any nef $a_1,\ldots,a_{d-1}$
as above, effectivity gives $\Theta_t\cdot\prod_j a_j\ge0$ by
approximation with ample classes. Since $\mu_*K_P=K_N$ as Weil cycles,
the projection formula gives
\begin{equation}\label{eq:subadjunction-big-decomposition}
 A_t\cdot\prod_j a_j
 \ge K_P\cdot\prod_j\mu^*a_j
 \ge c_S\prod_j(a_j^d)^{1/d}.
\end{equation}
Letting $t\downarrow0$ gives \eqref{eq:slc-stratum-mixed}. For the
top intersection take every $a_j=A_t$ in \eqref{eq:subadjunction-big-decomposition}.
Since $A_t$ is ample, division by $(A_t^d)^{(d-1)/d}$ gives
$(A_t^d)^{1/d}\ge c_S$. Letting $t\downarrow0$ now gives
$(D|_S)^d=B^d\ge c_S^d$, using $\deg\nu=1$.

Thus $D$ is nef and its restriction to every slc stratum is big, including
the zero-dimensional strata. This is log bigness for $[W,D]$, since its quasi-log
canonical strata are the slc strata and its non-quasi-log-canonical locus is empty
\cite[Theorem~2.9]{FujSub22}. For each slc stratum $S$ the open stratum is
$S\setminus\bigcup_{S'\subsetneq S}S'$, where $S'$ ranges over the smaller slc strata.
Hyperbolicity of $W$ makes every morphism from $\mathbb A^1$ to this open set
constant. Fujino's ampleness criterion gives the result
\cite[Definition~4.4 and Theorem~11.1]{FujMori25}.
\end{proof}

\section{Applications}\label{sec:conjecture-corollaries}

\subsection{Curvature conditions}\label{sec:curvature-corollaries}

For $(X_0,V_0):=(X,T_X)$, put $X_j:=\mathbb P(V_{j-1})$ and
$V_{j,(x,\ell)}:=(d\pi_j)^{-1}(\ell)$, where $\pi_j\colon X_j\to X_{j-1}$ and
$\ell\subset V_{j-1,x}$ is a line. Nondegenerate negative $k$-jet curvature means
that the tautological bundle $\mathcal O_{X_k}(-1)$ has a possibly singular
Hermitian metric $h_k$ whose dual has quasi-plurisubharmonic local weights and curvature satisfying $\Theta_{h_k^{-1}}[\xi]\ge\varepsilon|\xi|_{\omega_k}^2$
in the sense of currents for every smooth local section $\xi$ of $V_k$,
with some $\varepsilon>0$ and Hermitian form $\omega_k$ on $X_k$. The weights of $h_k$ must be locally bounded on the locus
of lifts of curve germs with nonzero first derivative \cite[Definitions~7.3, 7.4 and~7.12]{Dem97}.
\begin{corollary}\label[corollary]{cor:negative-directed-jets}
Let $X$ be a compact complex manifold with nondegenerate negative
$k$-jet curvature on $(X,T_X)$ for some $k\ge1$, in the sense of
\cite[Definition~7.12]{Dem97}. Then the conclusions of
\Cref{thm:main,cor:analytic-subspaces} hold for $X$.
\end{corollary}
\begin{proof}
Nondegeneracy and \cite[Theorem~7.8]{Dem97} give hyperbolicity, so
\Cref{thm:main,cor:analytic-subspaces} apply.
\end{proof}

\begin{corollary}\label[corollary]{cor:negative-finsler}
Let $G$ be a continuous complex Finsler metric on a compact complex
manifold $X$, smooth and positive off the zero section, with
$G(x,\lambda v)=|\lambda|^2G(x,v)$.
If its holomorphic sectional curvature $H_G$ is strictly negative in
every tangent direction, then the conclusions of \Cref{thm:main,cor:analytic-subspaces} hold for $X$.
\end{corollary}
\begin{proof}
If $X$ were not hyperbolic, Brody's lemma would give a nonconstant entire
curve $f$ with bounded derivative for a fixed Hermitian metric $h$
\cite{Brody78}, \cite[\S1, p.~1]{Duval21}. By homogeneity and compactness
of the $h$-unit sphere bundle, $u:=G(f,f')$ is bounded. Where $f'\ne0$,
the disc definition of $H_G$ gives
\begin{equation}
 -\frac{\Delta_{\rm eucl}\log u}{2u}=K_{f^*G}\le H_G(f,[f'])<0
\end{equation}
\cite[Definition~2.3 and (2.8)]{WanFins19}. Thus $u$ is strictly subharmonic
there and satisfies the submean inequality at its zeros since $u\ge0$.
This contradicts Liouville's theorem for bounded subharmonic functions on $\C$
\cite[proof of Proposition~3.1]{Dem97}. Apply \Cref{thm:main,cor:analytic-subspaces}.
\end{proof}

\subsection{Subvarieties and meromorphic maps}\label{sec:lang-projective}

\begin{corollary}\label[corollary]{cor:kob-lang}
Let $Y$ be a positive-dimensional irreducible closed subvariety of a
hyperbolic projective variety $W$, and let $\nu\colon Y^\nu\to Y$ be
its normalization. Then $Y$ is of general type. If $Y^\nu$ has log
canonical singularities, $K_{Y^\nu}$ is ample.
\end{corollary}
\begin{proof}
The projective variety $Y$ belongs to $\mathcal C$, so \Cref{thm:proj}
applies to the inclusion $Y\hookrightarrow W$. Its normalization is
projective since the finite map $\nu$ pulls back an ample bundle to an
ample bundle \cite[Proposition~1.2.13]{LazPosI}. Every entire curve
in $Y^\nu$ maps to a point in $W$ and hence lies in a finite fibre
of $\nu$. Thus $Y^\nu$ is hyperbolic, and in the log canonical case
\Cref{cor:slc-ample} applied to $(Y^\nu,0)$ gives ampleness.
\end{proof}

\begin{lemma}\label[lemma]{lem:general-type-image-lift}
Let $f\colon M\dashrightarrow W$ be a nonconstant meromorphic map from
a compact complex manifold to a compact hyperbolic complex space.
If $M\in\mathcal C$ or $W$ is projective, the reduced image $Y$ is of
general type and $f$ lifts to a dominant meromorphic map to a resolution of $Y$.

If $f\colon S\dashrightarrow W$ is a nonconstant rational map between
projective varieties with $S$ irreducible and $W$ hyperbolic, the
reduced image closure $Y$ is of general type and admits a projective
resolution to which $f$ has a dominant rational lift.
\end{lemma}
\begin{proof}
In the meromorphic case, the irreducible graph maps into one irreducible
component of $W_{\rm red}$, which is compact and hyperbolic. Kodama's theorem
for smooth sources makes $f$ holomorphic \cite[Theorem~1]{Kodama79}.
Its reduced image $Y$ is irreducible and closed analytic by Remmert's theorem
\cite[\S4, Satz~24]{Remmert57}, \cite[Ch.~II, Theorem~8.8]{DemCADG}.
If $M\in\mathcal C$, then $Y\in\mathcal C$ \cite[\S4.3, Lemma~4.6(1)]{Fujiki78}.
If $W$ is projective, Chow's theorem makes $Y$ projective \cite[Theorem~V]{Chow49}.
In either case apply \Cref{thm:proj} to $Y\hookrightarrow W$.
In the rational case the reduced image closure is projective, so \Cref{cor:kob-lang}
gives the same conclusion.

Choose a resolution $p\colon T\to Y$, projective in the rational case
\cite{Hir64}, and let $U\subset Y$ be the open set over which $p$ is an isomorphism.
The lift $p^{-1}\circ f$ is dominant and is defined where $f$ is defined and
takes values in $U$. In the holomorphic case its graph closure is the unique
irreducible component of $(M\times_Y T)_{\rm red}$ dominating $M$, proper over
$M$ and isomorphic over $f^{-1}(U)$. In the rational case $p^{-1}\circ f$ is rational.
\end{proof}
A compact manifold $M\in\mathcal C$ is \emph{weakly special} if no
connected finite \'etale cover of $M$ admits a dominant meromorphic map
to a positive-dimensional compact manifold of general type
\cite[Definition~9.26]{Cam04}.

\begin{corollary}
\label[corollary]{cor:weak-special-maps}\label[corollary]{cor:special}
Let $M'\to M$ be a connected finite \'etale cover of a positive-dimensional
weakly special compact manifold $M\in\mathcal C$. Every meromorphic map
$M'\dashrightarrow W$ to a compact hyperbolic complex space is constant.
\end{corollary}
\begin{proof}
A K\"ahler modification of $M$ pulls back to a modification of $M'$
whose source is K\"ahler, being finite \'etale over a K\"ahler manifold.
Thus $M'\in\mathcal C$. By \Cref{lem:general-type-image-lift}, a
nonconstant map $M'\dashrightarrow W$ would give a dominant meromorphic
map from $M'$ to a positive-dimensional resolution of general type,
contrary to weak specialness of $M$.
\end{proof}

A compact manifold $M\in\mathcal C$ is \emph{special} if every saturated
rank-one subsheaf $L\subset\Omega_M^p$, $p>0$, has $\kappa(M,L)<p$.
Here saturation means that $\Omega_M^p/L$ is torsion-free. A reduced irreducible closed analytic subspace of $M$
is special if a resolution of it is special \cite[Definition~2.24 and Theorem~2.27]{Cam04}.
Campana's core $c_M\colon M\dashrightarrow C(M)$ is the meromorphic fibration,
unique up to bimeromorphic equivalence, whose fibre through a very general
point is the largest such special subspace through that point
\cite[Theorem~3.3]{Cam04}. Here a very general point means a point outside
a countable union of proper closed analytic subsets. For an irreducible
projective variety $S$, Abramovich's Lang map
$\ell_S\colon S\dashrightarrow L(S)$ is the dominant rational map to a variety
of general type through which every dominant rational map from $S$ to a
variety of general type factors uniquely. Here $L(S)$ is a point if no such
map has positive-dimensional image \cite[Theorem~1.1 and Definition~1.2]{AbrLang}.
\begin{corollary}\label[corollary]{cor:special-maps}
Let $W$ be a compact hyperbolic complex space and $M$ a compact manifold
in class $\mathcal C$. Every meromorphic $f\colon M\dashrightarrow W$
factors through Campana's core $c_M$ by a unique meromorphic map
$C(M)\dashrightarrow W$. If $W$ is projective and $S$ is an irreducible
projective variety, every rational $f\colon S\dashrightarrow W$ factors
through Abramovich's Lang map $\ell_S$ by a unique rational map $L(S)\dashrightarrow W$.
\end{corollary}
\begin{proof}
A constant map factors through either base. For nonconstant $f$,
\Cref{lem:general-type-image-lift} gives a dominant lift $h$ to a
resolution $T$ of general type. In the meromorphic case $T\in\mathcal C$
by \cite[\S4.3, Lemma~4.6(3)]{Fujiki78}. Since $T$ is of general type, we may take $c_T=\mathrm{id}_T$. Then
functoriality of the core makes $h$ factor through $c_M$
\cite[Corollaries~3.11(1) and~3.14]{Cam04}. In the rational case take
a projective resolution $T$. The universal property of $\ell_S$ for
dominant rational maps to varieties of general type gives the
factorization of $h$ \cite[Theorem~1.1 and Definition~1.2]{AbrLang}.
Composing with $T\to Y\hookrightarrow W$, where $Y$ is the reduced
image closure, gives the desired factorization. Uniqueness follows
from dominance of $c_M$ and $\ell_S$.
\end{proof}

\section*{Acknowledgments}
The first author is grateful to Professor Fusheng Deng and Professor Xiaonan Ma
for their kind and thoughtful advice. The first author is supported by the Nankai Zhide Foundation,
the China Postdoctoral Science Foundation (CPSF, Grant Nos.~2025M783109
and~2026T190821), and the Postdoctoral Fellowship Program of CPSF
(Grade~B, Grant No.~GZB20260739).

\section*{AI use disclosure}
The authors developed the main ideas for the proofs of
\Cref{thm:main,thm:proj}. The disc argument for \Cref{thm:main}
builds on OpenAI's preprint \cite{OpenAI051} and on the authors'
earlier work \cite{QT26}.
They used ChatGPT to explore possible ways of weakening or removing the
hypotheses on currents in \cite[Theorem~A]{QT26}, namely $(C)$ for
nefness and $(C_2)$ for ampleness. ChatGPT also helped the authors check the proofs and simplify some of the arguments. 
The authors take full responsibility for the contents of this manuscript, including its mathematical arguments, references, and conclusions.

\end{document}